\documentclass[11pt,reqno]{amsart}
\usepackage{amssymb}
\usepackage{amsfonts}
\usepackage{amsthm}
\numberwithin{equation}{section}
\theoremstyle{plain}
\newtheorem{theorem}{Theorem}[section]
\newtheorem{proposition}[theorem]{Proposition}
\newtheorem{corollary}[theorem]{Corollary}
\newtheorem{lemma}[theorem]{Lemma}
\theoremstyle{definition}
\newtheorem{definition}[theorem]{Definition}
\newtheorem*{remark}{Remark}
\newtheorem*{example}{Example}
\newcommand{\refE}[1]{(\ref{E:#1})}
\newcommand{\refS}[1]{Section~\ref{S:#1}}

\newcommand{\refT}[1]{Theorem~\ref{T:#1}}

\newcommand{\refP}[1]{Proposition~\ref{P:#1}}
\newcommand{\refD}[1]{Definition~\ref{D:#1}}
\newcommand{\refC}[1]{Corollary~\ref{C:#1}}
\newcommand{\refL}[1]{Lemma~\ref{L:#1}}
\newcommand{\R}{\ensuremath{\mathbb{R}}}
\newcommand{\C}{\ensuremath{\mathbb{C}}}
\newcommand{\N}{\ensuremath{\mathbb{N}}}

\renewcommand{\P}{\ensuremath{\mathbb{P}}}
\newcommand{\Z}{\ensuremath{\mathbb{Z}}}
\newcommand{\K}{\ensuremath{\mathbb{K}}}

\newcommand{\cins}{\frac 1{2\pi\mathrm{i}}\int_{C_S}}

\newcommand{\g}{\mathfrak{g}}
\newcommand{\gb}{\overline{\mathfrak{g}}}
\newcommand{\gh}{\widehat{\mathfrak{g}}}
\newcommand{\Gb}{\overline{\mathcal{G}}}
\newcommand{\Gh}{\widehat{\mathcal{G}}}

\newcommand{\A}{\mathcal{A}}
\newcommand{\W}{\mathcal{W}}
\newcommand{\V}{\mathcal{V}}

\renewcommand{\L}{\mathcal{L}}

\newcommand{\res}{\operatorname{res}}

\renewcommand{\H}{\mathrm{H}}
\renewcommand{\a}{\ensuremath{\alpha}}

\newcommand{\Bh}{\widehat{B}}
\newcommand{\E}{\mathcal{E}}
\newcommand{\Dex}{D_1^*\cup D_{-1/2}^*\cup D_{-2}^*}

\newcommand{\tensor}{\otimes}   
\newcommand{\htensor}{\widehat{\otimes}}   
\begin{document}
\vspace*{-1cm}
\hbox{ }
{{\hspace*{\fill} math.QA/0412113}}

\vspace*{2cm}

\title[Global geometric deformations of current algebras] 
{Global geometric deformations of current algebras
as Krichever-Novikov type algebras}
\thanks{}
\author[A. Fialowski]{Alice Fialowski}
\address[Alice Fialowski]{Department of Applied Analysis,
E\"otv\"os Lor\'and University, P\'azm\'any P\'eter s\'et\'any 1,
H-1117 Budapest, Hungary}
\email{fialowsk@cs.elte.hu}
\author[M. Schlichenmaier]{Martin Schlichenmaier}
\address[Martin Schlichenmaier]{Universit\'e du Luxembourg,
Laboratoire de Mathematiques, Campus Limpertsberg, 
162 A, Avenue de la Faiencerie,
L-1511 Luxembourg}
\email{Martin.Schlichenmaier@uni.lu}
\begin{abstract}
We construct algebraic-geometric families of genus one (i.e. elliptic)
current and affine Lie algebras of Krichever-Novikov type.
These families deform the classical current, respectively 
affine Kac-Moody Lie algebras.
The construction is induced by the geometric process of
degenerating the elliptic curve to singular cubics.
If the finite-dimensional Lie algebra defining the
infinite dimensional current algebra is simple
then, even if restricted to local families, the constructed families are
non-equivalent to the trivial family.
In particular, we show that the current algebra is
geometrically not rigid, despite its 
formal rigidity.
This shows that in the infinite-dimensional Lie algebra case the
relations between geometric deformations,
formal deformations and Lie algebra two-cohomology are not that
close as in the finite-dimensional case.
The constructed families are e.g. of relevance in the global operator
approach to the Wess-Zumino-Witten-Novikov models 
appearing in the quantization of Conformal Field Theory.
\end{abstract}
\date{2.12.2004}
\maketitle

\vskip 1.0cm
\section{Introduction}\label{S:intro}

Deformations of mathematical structures 
is a tool of  fundamental importance in most parts of mathematics,
in mathematical physics, and physics.
Via deforming the object into a ``similar'' object, we get a richer
picture about the original object itself.
Moreover, via deformations we can approach the problem whether
we can equip the set of mathematical structures under consideration
(maybe up to certain equivalences) with the structure of
a topological or even geometric space.
In other words, does there exist a moduli space for these
structures?
If so, then for a fixed object the deformations of this object should
reflect the local structure of the moduli space at 
the point corresponding to this object.

This aspect of the  theory of deformations
 originated with the problem of classifying all possible pairwise
 non-isomorphic complex structures on a given differentiable
 real manifold. The fundamental idea, which should be credited to
 Riemann, was to introduce an analytic structure therein.

The notion of local and infinitesimal deformations of a complex analytic
 manifold first appeared in the work of Kodaira and Spencer.
In particular,
 they proved that infinitesimal deformations can be parameterized by a
corresponding
 cohomology group. The deformation theory of compact complex manifolds was
 developed further  by Kuranishi  and Palamodov.
 Shortly after the work of Kodaira and Spencer, algebro-geometric
 foundations were systematically developed by Grothendieck.
 For algebraic manifolds, the theory was constructed by Artin 
 and Schlessinger.

A typical example for results in these 
directions is the following. 
For the moduli space $\mathcal{M}_g$ of 
smooth projective curves of genus $g$ over $\C$ (or equivalently, compact
Riemann surfaces of genus $g$)
the tangent space $T_{[C]} \mathcal{M}_g$ can be naturally 
identified with 
$\H^1(C,T_C)$, where $T_C$ is the sheaf of holomorphic vector fields
over $C$.
This extends to higher dimension. In particular, it turns out that
for compact complex manifolds $M$,
the condition $\H^1(M,T_M)=\{0\}$ implies that $M$ 
is rigid, \cite[Thm. 4.4]{Kod}.
Rigidity means  that
any differentiable family $\pi:\mathcal{M}\to B\subseteq \R$, 
$0\in B$ which contains $M$ as the special member $M_0:=\pi^{-1}(0)$
 is trivial in a neighbourhood of $0$, i.e. for $t$ small enough
$M_t:=\pi^{-1}(t)\cong M$.
Even more generally, for a compact complex manifold $M$
and $\H^1(M,T_M)\ne \{0\}$ there exists a versal family 
which can be realized locally as a family  over a certain subspace
of    $\H^1(M,T_M)$ such that every appearing deformation family
is ``contained'' in this versal family
(see also \cite{Palm} for definitions,  results, and further references).

In this article we deal with  deformations of algebraic structures,
more specifically with  deformations of 
Lie algebras, in particular with such of infinite dimension.
 Formal deformations of arbitrary rings and
 associative algebras, and the related cohomology questions,
 were first investigated by Gerstenhaber, in a series of articles
 \cite{Ger}. The notion of deformation was applied to
 Lie algebras by Nijenhuis
 and Richardson \cite{NijRich1}, \cite{NijRich2}.

The cohomology space related to 
the deformations of a Lie algebra $\L$ is the Lie algebra
two-cohomology $\H^2(\L,\L)$ of $\L$ with  
values in the adjoint module.
Indeed, as long as this space is finite-dimensional,
it gives all infinitesimal deformations (i.e. there exists a universal
family of infinitesimal deformations).
Even more, in this case there exists a versal family for the 
formal deformations with base in this cohomology space,
see results of Fialowski \cite{Fiaproc}, Fialowski and Fuchs \cite{FiaFuMini},
from where \refT{formver} below is quoted.

The positive results in the deformation 
theory of compact analytic manifolds lead to the impression that the
vanishing of the relevant cohomology spaces will imply rigidity 
with respect to deformations also in the case of other
structures.
In an earlier paper \cite{FiaSchlV} the authors showed that this  hope is too
 naive.
There we considered the (infinite dimensional) Witt algebra (respectively its 
universal central extension, the Virasoro algebra) and constructed
algebraic-geometric deformations of it (each of them parameterized 
by the points of the affine line)
using Krichever-Novikov vector field algebras.
These deformations are non-trivial, only the special element 
in these families will
be isomorphic to the Witt algebra,
despite the fact that for the Witt algebra 
the two-cohomology space vanishes
\cite{Fiajmp}.
Hence, in the case of infinite dimensional Lie algebras 
the vanishing of the two-cohomology only  implies
infinitesimal  and formal rigidity.
Formal rigidity means that every deformation over the algebra
of formal power series in finitely many variables is equivalent to
a trivial family (see \refS{theory} for  precise definitions).
Our example shows that the Witt (or Virasoro) algebra despite
its formal rigidity is  geometrically not rigid.

This situation is peculiar for the case of infinite-dimensional
Lie algebras. For finite-dimensional Lie algebras we have 
strong relations between the different concepts of rigidity.
In particular, if a finite-dimensional 
Lie algebra is non-singular (see \refS{theory} 
for the definitions and further details)
it is infinitesimally, formally,  geometrically,
and analytically
rigid if and only if the two-cohomology space vanishes.

In this article we elaborate further on these phenomena.
We will consider the case of {\it current algebras}
$\gb=\g\otimes \C[z^{-1},z]$ and their central extensions
$\gh$, the {\it affine Lie algebras}.
Here $\g$ is a finite-dimensional Lie algebra.
Given an invariant, symmetric bilinear form $\beta$
the central extension $\gh$ is the vector space 
$\gb\oplus t\,\C$ endowed with the Lie bracket
\begin{equation*}
[x\otimes z^n, y\otimes z^m]=[x,y]\otimes z^{n+m}-\beta(x,y)\cdot 
n\cdot\delta_{m}^{-n}\cdot t,
\quad [t,\gh]=0,
\quad x,y\in\g,\ n,m\in\Z.
\end{equation*}
These algebras are of fundamental importance in a number of fields.
They supply examples of infinite dimensional algebras which are
still accessible to a structure theory. 
They appear as  gauge algebras in Conformal Field Theory
(CFT), \cite{BP}.
More generally, they are symmetry algebras of infinite-dimensional
systems, integrable systems, etc. 
More examples can be found in the book of Kac \cite{KacB}.
For  simple finite-dimensional Lie algebras $\g$ the central
extensions $\gh$, the associated affine Lie algebras,
are the {\it Kac-Moody algebras of untwisted affine type}.

For $\g$ finite-dimensional and simple
(and hence rigid) 
it was shown by Lecomte and Roger \cite{LeRoRid} that
the current algebra $\gb$ remains  formally rigid.
But we will exhibit again natural algebraic families of algebras 
containing $\gb$ as special element but all other members
will not be isomorphic
to $\gb$. In particular, these families can not be algebraic-geometrically 
equivalent to the trivial family.
Hence, the current algebra $\gb$ will  be geometrically 
not rigid.
All these families can be extended to families of centrally extended
algebras and yield in this way nontrivial deformations of
the affine algebra $\gh$.
The families constructed appear as families of higher-genus
multi-point current algebras of Krichever-Novikov type,
see \refS{kn} for their definition.
Hence, they are not just abstract families, but families
obtained by geometric process. 
The results obtained do not have only relevance to the
deformation theory of algebras, but are of importance in 
the fields in which current and affine algebras play a role.
Note that the maximal nilpotent Lie algebra of $\gb$ is even not
formally rigid
anymore, and its formal deformations are described in \cite{Fianil}. 

\medskip 
In particular, they are of 
relevance in 
two-dimensional CFT and its quantization.
It is well-known that the Witt algebra, 
the Virasoro algebra, and their representations 
are of fundamental importance  for 
conformal field theory on the Riemann sphere
(i.e. for genus zero), see \cite{BP}.
Krichever and Novikov \cite{KNFac} proposed in  case of
higher genus Riemann surfaces (with two insertion 
points) the use of global operator fields 
which are given with the help of  the Lie algebra 
of vector fields of Krichever-Novikov type,
certain related algebras, 
and their representations (\refS{kn}).

Their approach was extended to the multi-point
situation (i.e. an arbitrary number of insertion points was allowed) by
Schlichenmaier
\cite{SchlDiss}, \cite{Schlkn}, \cite{Schleg}, \cite{Schlce}.
The necessary central extensions where constructed and higher
genus multi-point current and affine algebras were introduced
\cite{SchlCt}.
These algebras consist of meromorphic objects on a Riemann surface
which are holomorphic outside a finite set $A$  of points. In
turn the set $A$ is divided into two disjoint subsets
$I$ and $O$. With respect to some possible interpretation of the
Riemann surface as the world-sheet of a string, the points
in $I$ are called {\it in-points}, the points in $O$ 
are called {\it out-points}, 
corresponding to incoming and outgoing free strings;
the world-sheet itself corresponds to possible interaction.
This splitting introduces an almost-graded structure (see \refS{kn})
for the algebras and their representations.
Such an almost-graded structure is needed to construct representations
of relevance in the context of CFT, e.g. highest weight
representations, 
fermionic Fock space representations, etc.
  
In the process of quantization of  conformal fields one has
to consider families of algebras and representations over the
moduli space of compact Riemann surfaces
(or equivalently, of  smooth projective curves over $\C$) of genus $g$ with
$N$ marked points.
Models of most importance in CFT are the 
Wess-Zumino-Witten-Novikov models (WZWN).
Tsuchiya, Ueno, and Yamada \cite{TUY} gave a sheaf version of
WZWN models over the moduli space.
In \cite{SchlShWz}, 
\cite{SchlShWz2}
Schlichenmaier and Sheinman  developed a global operator version.
In this context
of particular interest is the situation
$I=\{P_1,\ldots, P_K\}$, the marked points we want to vary,
 and $O=\{P_\infty\}$, a reference point. We obtain
families of algebras over 
the moduli space $\mathcal{M}_{g,K+1}$ of curves of genus $g$ with
$K+1$ marked points, and we are exactly in the middle of the
main subject of this article.
In \cite{SchlShWz} and   \cite{SchlShWz2} it is shown that
there exists a global operator description of 
WZWN  models with the help of the Krichever
Novikov objects at least over  a dense open subset of the moduli space.
The following is  just a very rough outline.
Let us start from families of representations $\V$ of families of
higher genus affine algebras 
(see \refS{kn} for their definition). 
The vector bundle  of conformal blocks can 
be defined as the vector bundle with fiber 
(over the moduli point $b=[(M,\{P_1,\ldots, P_K \},\{P_\infty\})]$)
as follows. Take the quotient of  the fiber of the representation
$\V_b$ by  the subspace generated by
the vectors obtained by the action
 of those elements of the affine algebra  which vanish at the
reference point $P_\infty$ (i.e. the fiber of the quotient is the space of 
coinvariants of this subalgebra).

The bundle of conformal blocks carries a connection
called
the {\it Kniz\-hnik-Zamolod\-chikov connection}.
In its definition the Sugawara construction plays
an important role, which 
associates to representations of affine algebras representations
of the (almost-graded) centrally extended vector
field algebras, see \cite{SchlShSug}.
A certain subspace of the vector field algebra  (assigned
to the moduli point $b$) corresponds to tangent directions
on the moduli space $\mathcal{M}_{g,K+1}$ at the point $b$.

Now clearly, the following question is fundamental. What happens
if we approach the boundary of the moduli space?
The boundary components correspond to curves with singularities.
Resolving the singularities yields curves of lower genera.
By geometric degeneration we obtain families of (Lie) algebras containing
a lower genus algebra (or sometimes a subalgebra of it), 
corresponding to a suitable collection of marked points,
as special element.
Or reverting the perspective,
we obtain a typical situation of the deformation  of an algebra
corresponding in some way to a lower genus situation, containing
higher genus algebras as the other elements in the family.
Such kind of geometric degenerations
are fundamental if one wants to prove Verlinde
type formula via factorization and normalization technique,
see \cite{TUY}.

By a maximal degeneration a collection of $\P^1(\C)$'s will 
appear.
Indeed, the examples considered in this article are exactly of this type.
The deformations  appear as families 
of current algebras which are naturally defined over the
moduli space 
of genus one curves (i.e. of elliptic curves, or
equivalently of complex one-dimensional tori) with two marked points.
These deformations are
associated to  geometric degenerations of elliptic curves to
singular cubic curves. 
The desingularization (or normalization) of their singularities will yield 
the projective line as normalization. We will end up with
algebras related to the genus zero case.
The full geometric picture behind the degeneration 
 was discussed in \cite{SchlDeg}.
In particular, we want to point out, that even if one starts with
two marked points, by going to the boundary of the moduli space
one is forced to consider more points (now for a curve of lower
genus).

In special cases 
the classical current (or affine) algebras 
appear as 
degenerations of elliptic two-point current algebras. 
Considered from the opposite point of view,
in the sense of this article,  
the elliptic two-point current algebras
are global deformations of the classical current algebra.
Nevertheless, as we show here, the structure of 
these algebras are not determined 
by the classical current  algebras, despite the formal rigidity of
the latter (if $\g$ is simple), \cite{LeRoRid}.

\medskip

The structure of the article is the following.
In \refS{theory} we recall the different definitions  of 
deformations of Lie algebras.
In particular we stress the fact, that 
in the case of infinite dimensional Lie algebra it is important
 to distinguish clearly the different notions, and
 indicate in which category we work:
infinitesimal, formal, geometric, or analytic deformations.
Also it is important to allow as base of the deformations not
only local algebras, but also global algebras.
They correspond to geometric situations. The most simple global case is
the case of the algebra of polynomials in one variable.
Geometrically this corresponds to  a deformation over the affine line.
As already indicated above, in the formal case everything can be 
described in cohomological terms.
To contrast the infinite dimensional case with the finite-dimensional
one, we recall some results from their theory.
In the finite-dimensional case for non-singular Lie algebras
all notions of rigidity introduced above are equivalent and
correspond to the fact that the cohomology space vanishes.
The reason is that in this case the whole situation can 
be described within the frame of finite-dimensional algebraic
geometry.

\medskip
In \refS{kn} we recall 
what is needed about the higher genus multi-point algebras of
Krichever-Novikov type.
The following algebras are introduced: the associative algebra
of functions and the Lie algebras of vector fields and of currents,
including
their central extensions. For the currents the central extensions are the
higher genus multi-point affine algebras.

\medskip
In \refS{elliptic} 
we construct  geometric deformations of the standard current
algebra 
by considering  certain families of algebras for
the genus one case (i.e. the elliptic curve case) 
and let the elliptic curve degenerate to
a singular cubic.
The two points, where poles are allowed, are the zero element of the
elliptic curve (with respect to its group  structure) and a
2-torsion point.
In this way we obtain families parameterized over the affine line with
the peculiar behaviour that every family is a global deformation of the
classical current  algebra, 
i.e. the classical current algebra  
is a special member, whereas all other members
are mutually 
isomorphic but not isomorphic to the special element if the 
finite-dimensional Lie algebra is simple,
 see \refT{families}. Even if
restricted to  small open neighbourhoods of the point
corresponding to the special element, 
these families are non-trivial, only  infinitesimally and formally they
are trivial. 
The construction can be extended to the centrally extended 
algebras, yielding  global deformations of 
the affine algebra.

\medskip
In \refS{degen} we consider the geometric picture behind it.
In particular, we identify those algebras we obtain over the nodal 
cubics (i.e. the cubic curves with one singular point with two 
tangent directions at this point).
We explain the geometric reason why we obtain them.
Depending whether the node will become the point where a pole is
allowed or not, we obtain a three-point current algebra of genus zero or
a certain subalgebra of the classical current algebra.

\medskip
In   \refS{cohomo} we give
the cohomology classes of our families of deformations.
As it is known that these algebras are formally rigid in the 
simple case, it can be expected that, in general, the
cocycles will be coboundaries. We show this by direct calculations.

\medskip
In an appendix we calculate the cocycle defining the central 
extensions of our family.

\medskip
In the process of constructing our families we constructed families
of the commutative and associative algebra of functions on
Riemann surfaces with prescribed regularity.  
These are families deforming the algebra of Laurent polynomials
$\C[z^{-1},z]$.
This algebra is the coordinate algebra of the smooth affine curve
 $\P^1\setminus\{0,\infty\}$.
The  cohomology 
corresponding to commutative and associative deformations
(i.e. the Harrison cohomology)
vanishes for such algebras. The algebra 
is infinitesimally
and formally rigid. Nevertheless, we constructed
nontrivial (geometrically) local  families.
Maybe, Kontsevich's concept of semi-formal
deformations \cite{Kon} (related to filtrations
of certain type) might help to 
describe this situation better.
Indeed, by the almost-grading of the  families 
considered in this article, the families  of associative 
algebras are semi-formal
deformations in his sense.

\medskip
{\it Acknowledgements:}
The authors thank different institutions for hospitality experienced 
during the preparation of this article.
A.~F. and M.~Sch. thank the Erwin-Schr\"odinger Institute (ESI) in Vienna,
M.~Schl. the E\"otv\"os Lor\'and University in Budapest,  and the
Institut des Hautes \'Etudes Scientifiques (IHES)
in Bures-sur-Yvette.
Discussions with M. Kontsevich, P. Michor and E. Vinberg 
are gratefully acknowledged.
The work was partially supported by
grants OTKA T034641 and T043034.


\section{Some generalities on deformations of Lie algebras}\label{S:theory}
\subsection{Intuitive description}\label{S:intuitive}
$ $

Let us start with the intuitive definition
of a Lie structure depending on a parameter $t$. 
Let $\L$ be a Lie algebra
with Lie bracket $\mu=\mu_0$ over a field $\K$. A {\it deformation} of
$\L$ is a one-parameter family $\L_t$ of Lie algebras 
(with the same underlying vector space) with the bracket
\begin{equation}\label{E:intdef}
\mu_t = \mu_0 + t\phi_1 + t^2\phi_2 + ...
\end{equation}
where $\phi_i$ are $\L$-valued (alternating) two-cochains, i.e. elements of
$\mathrm{Hom}_{\K}(\bigwedge^2\L,\L)=C^2(\L,\L)$, and $\L_t$ is a Lie
algebra for each $t\in \K$. (see \cite {Fiasbor, Ger}). 
In particular, we have $\L=\L_0$. Two deformations
$\L_t$ and $\L'_t$ are equivalent if there exists a linear automorphism
$\hat\psi_t = \text{id} + \psi_1t + \psi_2t^2 + ...$ of $\L$ where $\psi_i$
are linear maps over $\K$, i.e. elements of $C^1(\L,\L)$, such that
\begin{equation}\label{E:intequ}
\mu'_t(x,y) = \hat\psi_t^{-1}(\mu_t(\hat\psi_t(x),\hat\psi_t(y))).
\end{equation}

These objects are related to  Lie algebra cohomology.
We recall the following 
definitions.
A bilinear map $\omega:\L\otimes \L\to\L$ is a 
Lie algebra two-cocycle with values in the adjoint representation
$\L$, if $\omega$  is alternating  and fulfills
\begin{multline}
0=d_2(\omega)(x,y,z):=
\omega([x,y],z)-\omega([x,z],y)+\omega([y,z],x)
\\
-[x,\omega(y,z)]+[y,\omega(x,z)]-[z,\omega(x,y)].\qquad
\end{multline}
The vector space  of two-cocycles is denoted by
$\mathrm{Z}^2(\L,\L)$.
A two-cocycle $\omega$ is a coboundary if
there exists  a linear map $\eta:\L\to\L$
such that
\begin{equation}
\omega(x,y)=(d_1\eta)(x,y):=
\eta([x,y])-[x,\eta(y)]-[\eta(x),y].
\end{equation}
The vector space of coboundaries is 
denoted by
$\mathrm{B}^2(\L,\L)$ and is a  subspace of $\mathrm{Z}^2(\L,\L)$.
The quotient space is the Lie algebra two-cohomology
with values in the adjoint module, 
denoted by  $\H^2(\L,\L)$, see e.g. \cite{Fuks}, \cite{FeFu} for
further details.

Coming back to the deformation \refE{intdef} we see that
the Jacobi identity for the algebras $\L_t$ implies that the two-cochain
$\phi_1$ is indeed a cocycle, i.e. it fulfills
$d_2\phi_1=0$.
If $\phi_1$ vanishes identically, the first nonvanishing $\phi_i$ will
be a cocycle.
If $\mu'_t$ is an equivalent deformation (with cochains $\phi'_i$) then
\begin{equation}
\phi'_1-\phi_1=d_1\psi_1.
\end{equation}
Hence, every equivalence class 
of deformations defines uniquely an element of $\H^2(\L,\L)$.
This class is called the {\it differential} of the deformation.
The differential of a family which is
equivalent to a trivial family will be the zero cohomology class.
\subsection{Global deformations} 
$ $

The deformation $\L_t$ might be considered not only 
as a family of Lie algebras, but also as a Lie algebra over the 
algebra $\K[[t]]$ of formal power series over 
$\K$. 
From this description a natural step is to allow  
$\K[[t_1,t_2,\ldots, t_k]]$,  or 
an arbitrary  commutative algebra $A$ 
over $\K$ with unit  as the base
of a deformation. 

In the following we will assume that $A$ is a
commutative algebra over  $\K$ (where $\K$ is
a field  of  characteristic zero)
which  admits an augmentation 
 $\epsilon: A \to \K$.
This says that $\epsilon$ is a $\K$-algebra homomorphism, e.g.
$\epsilon(1_A)=1$.
The ideal $m_\epsilon:= \text{Ker\,} \epsilon$ is 
a maximal ideal of $A$. 
Vice versa, given a maximal ideal $m$ of $A$ with 
$A/m\cong\K$, the natural quotient map defines an
augmentation. 

If $A$ is a finitely generated $\K$-algebra over
an algebraically closed field $\K$ then $A/m\cong\K$
is true for every maximal ideal $m$. Hence, in this case 
every such $A$ admits at least one augmentation and
all maximal ideals are coming from augmentations.

Let us consider a Lie algebra $\L$ over the field $\K$, 
 $\epsilon$ a fixed augmentation of $A$, 
and
$m=\text{Ker\,} \epsilon$ the associated maximal ideal. 
\begin{definition}
\label{D:glob} (\cite{FiaFuMini})
 A {\it global deformation} $\lambda$ of $\L$ with 
base $(A,m)$ or simply with  base $A$, is a Lie $A$-algebra structure
on the tensor product $A\tensor_\K\L$ with bracket $[.,.]_\lambda$ such
that
\begin{equation}
\epsilon \textstyle{\tensor} \text{id} : A \textstyle{\tensor} \L \to
 \K \textstyle{\tensor} \L = \L
\end{equation}
is a Lie algebra homomorphism.

Specifically, it means that for all $a,b \in A$ and $x,y \in \L$,
\smallskip

\quad (1) $[a \tensor x, b \tensor y]_\lambda = (ab \tensor \text{id\,})[1 \tensor x,
1 \tensor y]_\lambda$,
\smallskip

\quad (2) $[.,.]_\lambda$ is skew-symmetric and satisfies the Jacobi 
identity,
\smallskip

\quad (3) $\epsilon \tensor \text{id\,}([1 \tensor x, 1 \tensor y]_\lambda) = 1
\tensor [x,y]$.
\end{definition}

\noindent 
By Condition (1) to describe a deformation 
it is enough to give  the elements $[1 \tensor x, 1
\tensor y]_\lambda$ for all $x,y \in \L$. 
From condition (3) it follows  that for them the Lie product has the form
\begin{equation}\label{E:defex}
[1 \textstyle{\tensor} x, 1 \textstyle{\tensor} y]_\lambda = 1
\textstyle{\tensor} [x,y] + \sum_i a_i \textstyle{\tensor} z_i,
\end{equation}
with $a_i \in m$, $z_i \in \L$.

A deformation is called {\it trivial} if $A\tensor_{\K} \L$
carries the trivially extended Lie structure, i.e. \refE{defex} reads as 
$[1\tensor x,1\tensor y]_\lambda=1\tensor[x,y]$.
Two deformations of a Lie algebra $\L$ with the same base $A$ are
called {\it equivalent} if there exists a Lie algebra isomorphism
between the two copies of $A \tensor \L$ with the two Lie algebra
structures, compatible with $\epsilon \tensor \text{id}$.

We say that  a deformation
is {\it local (in the algebraic sense)} 
if $A$ is a local $\K$-algebra 
with unique maximal ideal $m_A$. 
By assumption $m_A=\text{Ker\,} \epsilon$
and $A/m_A \cong \K$.
In case that in addition  ${m_A}^2 = 0$, the deformation
is called  {\it infinitesimal}.
\medskip

An important class of examples is given if $A$ is the algebra
of regular functions of an affine variety $V$, i.e.
$A=\K[V]$. The algebra $A$ is also called the coordinate algebra of
the affine variety $V$.  
Let us assume that $\K$ is algebraically closed.
It is known that every finitely generated (as ring) 
$\K$-algebra which is reduced, i.e. which has no nilpotent elements,
is the algebra of regular functions of a suitable 
affine variety. The variety might be reducible if $A$ has zero-divisors.
The  maximal ideals $m_x$ of $A$ correspond exactly to
the points $x\in X$, and exactly to the set of augmentations
$\epsilon_x$.
Fixing an augmentation $\epsilon_{x_0}$ means fixing a point $x_0\in V$.
If $A$ is a non-local ring, there will be different maximal ideals,
hence also different augmentations.
Let $\L$ be a $\K$-vector space and assume that there exists a 
Lie $A$-algebra structure $[.,.]_A$ on $A\tensor_{\K}\L$.
Given an augmentation $\epsilon:A\to\K$ with 
associated maximal ideal  $m_\epsilon=\text{Ker\,} \epsilon$,
one obtains a Lie $\K$-algebra structure $\L^\epsilon=(\L,[.,.]_\epsilon)$
 on the
vector space $\L$ by passing to the quotient
$A/m_\epsilon$.
In this way the Lie algebra $A\tensor_{\K}\L$ gives a family
of Lie algebra structures on the space $\L$ parameterized by the points
of $V$. We obtain a map
\begin{equation}
V\quad \to\quad 
\{\text{ set of Lie algebra structures on the vector space $\L$ } \}.
\end{equation}

We will denote both $A=\K[V]$ and $V$ as the base of the
deformation. It is quite convenient to consider both  
(the algebraic and the geometric) pictures.
Clearly, the Lie $A$-algebra $A\tensor_{\K}\L$, i.e. the
global family over $A$, is a 
deformation in the sense of \refD{glob}
with basis $(A,m_x)$ for all points $x\in V$.

If we choose a basis $\{T_a\}_{a\in J}$ 
of the Lie algebra $\L$, the Lie structure in $\L$ is given
by the structure constants $\{C_{a,b}^c\}$ defined via
\begin{equation}
[T_a,T_b]=\sum_{c\in J}{}' C_{a,b}^cT_c,\qquad  a,b\in J.
\end{equation}
Here $\sum'$ denotes, that only a finite number of the summands
will be different from 0.
Assume that we have a deformation over $A=\K[V]$ in the sense
of \refD{glob}, where $\L$ lies above $x_0$.
The elements $1\otimes T_a$ are now an $A$-basis of
$A\otimes \L$ and Equation \refE{defex} can be written
as
\begin{equation}
[1\otimes T_a,1\otimes T_b]_{\epsilon_x}=
\sum_{c\in J}{}' C_{a,b}^c(x)\, 1\otimes T_c\qquad  a,b\in J
\end{equation}
with algebraic functions
\begin{equation}
C_{a,b}^c(x)\in\K[V], \quad
\text{with}\quad
(C_{a,b}^c(x)-C_{a,b}^c)\in m_{x_0},
\text{ or equivalently }
C_{a,b}^c(x_0)=C_{a,b}^c.
\end{equation}
Note that the range of the summation might be bigger than
the original one.
In this way we get back the intuitive picture of algebraically varying
structure constants. The structure constants of the Lie algebra 
$\L_{x'}$ lying
above the point $x'$ are given by  $C_{a,b}^c(x')$.

A very import special case is a
deformation over the affine line $\mathbb{A}^1$.
 Here the corresponding algebra is $A=\K[\mathbb{A}^1]=\K[t]$, the
algebra of polynomials in one variable.
For a deformation of the Lie algebra $\L=\L_0$ over the affine
line, the Lie structure $\L_\alpha$ in the fiber over the point $\alpha\in\K$ 
is given by considering the augmentation corresponding to the maximal
ideal $m_\alpha=(t-\alpha)$. 
\begin{remark}
Despite the fact that we consider only deformations over affine
varieties, 
we call this type of deformations {\it global}, as we 
allow also non-local algebras.
Nevertheless, it might be useful to study deformations over more
general global basis, like projective varieties, schemes, complex
manifolds, analytic spaces, etc. 
As in this article we  are only interested in rigidity 
questions, and they are geometrically local, 
i.e. only the case of open neighbourhoods (e.g. local affine neighbourhoods)
is of importance, see \refD{rigidalg},
the given definition 
will be sufficient for us.
\end{remark}

\medskip

For further reference we note
\begin{definition}\label{D:push}
Let $A'$ be another commutative algebra over $\K$ with a fixed
augmentation $\epsilon': A' \to \K$, and let $\phi: A \to A'$ be an
algebra homomorphism with $\phi(1) = 1$ and $\epsilon'\circ\phi =
\epsilon$. If a deformation $\lambda$ of $\L$ with base
$(A,\text{Ker\,}\epsilon = m)$ is given, then the {\it push-out}
$\lambda'=\phi_*\lambda$ is the deformation of $\L$ with base $(A',
\text{Ker\,}\epsilon' = m')$, and Lie algebra structure
$$
{[a_1' \textstyle{\tensor}_A (a_1 \textstyle{\tensor} l_1),
a_2' \textstyle{\tensor}_A (a_2 \textstyle{\tensor} l_2)]}_{\lambda'}
 := a'_1a'_2
\textstyle{\tensor}_A {[a_1 \textstyle{\tensor}l_1,a_2 \textstyle{\tensor}
l_2]}_{\lambda},
$$
($a'_1, a'_2 \in A', a_1, a_2 \in A, l_1, l_2 \in \L$)
on $A' \tensor \L = (A' \tensor_A A) \tensor \L = A' \tensor_A (A \tensor
\L)$.  Here $A'$ is regarded as an $A$-module with the structure $aa' =
\phi(a)a'$.
\end{definition}
\subsection{Formal deformations}
$ $

Let $A$ be a complete local algebra over $\K$, so
$A = \overleftarrow{\lim\limits_{n\to \infty}}(A/m^n)$, where $m$ is
the maximal ideal of $A$. Furthermore,  we 
will assume that $A/m\cong\K$, and $\dim (m^k/m^{k+1})<\infty$ for all $k$.
\begin{definition}
A {\it formal deformation} of $\L$ with base $A$ is a Lie algebra
structure on the completed tensor product
$A \htensor\L =
    \overleftarrow{\lim\limits_{n\to \infty}}((A/m^n)\tensor \L)$
such that
\begin{equation}
\epsilon\textstyle{\htensor}
\text{id}: A \textstyle{\htensor}\L \to \K\tensor \L = \L
\end{equation}
is a Lie algebra homomorphism.
\end{definition}

If $A=\K[[t]]$, then a formal deformation of $\L$ with base $A$
is the same as a formal 1-parameter deformation of $\L$ (see \cite{Ger}).

There is an analogous definition for equivalence of deformations
parameterized by a complete local algebra. 
\subsection{The case of finite-dimensional Lie algebras}\label{S:finite}
$ $

In this subsection let $\L$ be a finite-dimensional Lie algebra of
dimension $n$ over $\K$.
Let $\{T_a\}_{a=1,\ldots, n}$ be a set of 
basis elements for the vector space $V$ the Lie algebra $\L$ is modeled on.
The Lie structure $\mu$ is fixed by  the structure constants
$\{C_{a,b}^c\}_{a,b,c=1,\ldots, n}$,
defined by
\begin{equation}\label{E:lieb}
[T_a,T_b]=\sum_{c=1}^n C_{a,b}^c T_c,\qquad  a,b=1,\ldots,n.
\end{equation}
From the skew-symmetry and the Jacobi identity for the Lie algebra
$\L$ it follows that the structure constants obey the following
relations
\begin{equation}\label{E:structlie}
\begin{gathered}
C_{a,b}^c+C_{b,a}^c=0,\quad  a,b,c=1,\ldots,n,
\\
\sum_{l=1}^n
\left(
C_{a,b}^l C_{l,c}^d+
C_{b,c}^l C_{l,a}^d+
C_{c,a}^l C_{l,b}^d
\right)=0,
\quad  a,b,c,d=1,\ldots,n.
\end{gathered}
\end{equation}
Conversely, given such a set of $\{C_{a,b}^c\}_{a,b,c=1,\ldots, n}$,
satisfying \refE{structlie}
the vector space generated by $n$ basis elements $T_a$ carries
the structure of a Lie algebra defined via \refE{lieb}.

Obviously, \refE{structlie} are  algebraic equations and
the vanishing set $Lalg_n$ 
in $\K^{N}$, $N=n^3$, of the ideal generated by them 
``parameterizes'' the possible Lie algebra structures on
the $n$-dimensional vector space $V$.
Strictly speaking, it is more appropriate to talk about the
scheme, as one should better consider the not necessarily
reduced structure on  $Lalg_n$, see Rauch \cite{Rau}.
Furthermore, as $\mu$ is a bilinear map
$V\times V\to V$,
the structure constants might be considered more
canonically as elements of $V^*\otimes V^*\otimes V$ with
 $V^*$ denoting the dual space of $V$.

Given a linear automorphism  $\Phi\in \mathrm{GL}(V)$,
it will define 
an action   on  $V^*\otimes V^*\otimes V$ by
\begin{equation}
(\Phi\star\mu)(x,y)=\Phi(\mu(\Phi^{-1}(x),\Phi^{-1}(y))),
\end{equation}
which respects  the conditions defining a Lie algebra 
structure.
The Lie algebras $(V,\mu)$ and  $(V,\mu')$
will be isomorphic 
iff  $\mu$ and $\mu'$ are  
in the same orbit under the  $\mathrm{GL}(V)$ action.
On the level of structure constants, i.e. after choosing a
basis in $V$, we obtain a $\mathrm{GL}(n)$ action on
 $Lalg_n$.
In this way the isomorphy classes of Lie algebras of dimension
$n$ correspond exactly to the  $\mathrm{GL}(n)$ orbits of $Lalg_n$.

For further reference we note, that a finite-dimensional Lie algebra $\L$ is 
called nonsingular if the corresponding point $\mu$ in  
the scheme $Lalg_n$ is 
a non-singular point in the sense of algebraic geometry.
\subsection{Rigidity}
$ $

Intuitively, rigidity of a Lie algebra  $\L$ means that we 
cannot deform the Lie algebra. Or, formulated differently, given a family
of Lie algebras containing $\L$ as the special element $\L_0$, any 
element $\L_t$ in the family ``nearby'' will be isomorphic to
  $\L_0$.
Of course, the definition depends on the category in which the
deformations are considered.

As ``nearby'' is already encoded in the infinitesimal and formal 
setting we can define:
\begin{definition}\label{D:rigidform}
(a) A Lie algebra $\L$ is {\it infinitesimally rigid} iff 
every infinitesimal deformation of it is equivalent to  the trivial 
deformation.
\newline
(b) A Lie algebra $L$ is {\it formally rigid} iff 
every formal deformation of it is equivalent to  the trivial deformation.
\end{definition}
For the global deformation such a condition  would be too much 
to require as 
the following example shows.
\begin{example}
Let $\L$ be any non-abelian Lie algebra with Lie bracket $[.,.]$.
We define a family of Lie algebras $\L_t$ over 
the algebra $\K[t]$ 
(i.e. geometrically over the affine line $\K$) by
taking as Lie bracket the bracket
$[x,y]_t:=(1-t)[x,y]$.
If we set $t=0$ we obtain back our Lie algebra $\L$. 
Moreover, as long as $t\ne 1$ the algebras $\L_t$ are isomorphic to
$\L$, but $\L_1$ as abelian Lie algebra will be non-isomorphic.
Hence such a family will never be the trivial family.
But if we restrict the family to  the 
(Zariski) open subset $\K\setminus\{1\}$,  of the base $\K$, 
we obtain a trivial 
family.
\end{example}

In the purely algebraic setting we can apply Zariski topology.
The Zariski topology (Z-topology) of a ring $R$ is  a topology on 
the set of prime ideals, yielding a topological space $Spec(R)$.
Instead of considering the general case, we
will restrict ourselves to the 
situation  which is of interest here.
Let $\K$ be an algebraically closed field of characteristic zero and
$A$ a finitely generated reduced $\K$-algebra.
Recall that this corresponds geometrically to the case of
 affine  varieties $V$ as base. 
A subset $W$ of $V$ is called Z-closed if it is the vanishing
set of an ideal in $A=\K[V]$. The set $W$ itself will be an affine
variety.
A subset $U$ of $V$  is called Z-open if it is the
complement of a closed set.
A basis of the open set of the topology is given by Z-open
affine varieties. These basis sets are open sets in $V$ 
and they themselves are affine (this says  they are the vanishing sets
of an ideal in a suitable affine space), see \cite{Harts} for
this and further results.

In particular, given an arbitrary Z-open subset $U$ containing the
point $x\in V$ we can always find an open subset $W$ of $U$ which 
itself is affine.
For simplicity (and in accordance with the  general picture)
let us use $Spec(A)$ to denote the variety $V$ and let us identify
the (closed) points in $V$ with the maximal  ideals.

\begin{definition}\label{D:rigidalg}
A Lie algebra $\L$ is {\it algebraic-geometrically rigid} 
of just {\it geometrically rigid} if for every
deformation of $\L$ in the sense of \refD{glob}
with base $(A,m)$, where $A$ is a 
finitely generated reduced 
 algebra 
over the algebraically closed field $\K$ and $m$ a
maximal ideal, there exists a Zariski open affine
neighbourhood $Spec(B)$  of the point $m$ in $Spec(A)$ 
such that 
the restricted family over $(B,m)$ is equivalent to the trivial 
family.
\end{definition}
Note that the restriction  given by the geometric picture
$Spec (B)\hookrightarrow Spec (A)$ corresponds exactly to the 
push-out of the deformation given by the induced
$\K$-algebra homomorphism $A\to B$.
Note also, that with a slight abuse of notation we understand
in  $(B,m)$ by $m$ the ideal generated by the image in $B$
of the ideal $m$ of $A$.

In the case if the base field is $\C$ (or $\R$)
and the base of the deformation is a finite-dimensional
analytic manifold $M$ with a chosen point $x_0\in M$
we consider also the usual topology on $M$.
In particular, we can talk about open subsets of $M$ containing
$x_0$.
\begin{definition}\label{D:rigidana}
A Lie algebra $\L$  over $\C$ (or $\R$) is {\it analytically rigid} if 
for every family over a finite-dimensional analytic manifold $(M,x_0)$ 
with base point $x_0$,  with  special fiber  $\L\cong \L_{x_0}$ 
over the point 
$x_0\in M$, there is an open neighbourhood $U$ of $x_0$, such that
the restriction of this family is equivalent to the  trivial family.
\end{definition}
Obviously, as the analytic rigidity is a (geometrically) local
condition,
it is enough to establish rigidity by considering families over
$\C^k$ (resp. $\R^k$).

There are other definitions of rigidity by considering
differentiable manifolds or analytic spaces as base.
Moreover, we could incorporate infinite dimensional manifolds or 
non-Noetherian algebras as base of the deformations.
But  for the families discussed in
this article we do not need them. Hence, we will not discuss them here.
Already in the case that the Lie algebra is infinite dimensional and
the base is finite-dimensional, interesting examples appear.
It can be expected that the situation will become even more
intricate by allowing infinite dimensional bases.
To our knowledge not much is known in this direction, whereas
the topic is surely interesting. We intend to come back to it
in the future.

\bigskip

In case that the Lie algebra is {\bf finite-dimensional} there is another
important definition of rigidity.
\begin{definition}\label{D:rigidorbit}
(\cite{Ger}, \cite{NijRich2})
Let $\L$ be a  Lie algebra of dimension $n$ corresponding to
the  point $\mu=\{C_{a,b}^c\}$ in the space of structure constants.
The algebra $\L$ is called {\it rigid (in the orbit sense)} or
just {\it rigid} if the orbit of  $\{C_{a,b}^c\}$ under
the $\mathrm{GL}(n)$ action is Zariski open in $Lalg_n$.
\end{definition}
Given a finite-dimensional Lie algebra $\L$ which is  rigid 
(in the orbit sense) and
 a deformation  over the base $(V,x_0)$, then the algebras $\L_x$ will
 be isomorphic to $\L_{x_0}=\L$ for any x in a Zariski open
 neighbouhood $U$ of $x_0$.
Hence, (orbitally)  rigidity implies 
geometric (and if $\K$ is equal $\R$ or $\C$ also analytic) rigidity.
See \refT{fd0} and \refT{fdrig} for partial inversions of this statement.

\subsection{Rigidity and cohomology}
$ $

\label{S:H2}
As explained in \refS{intuitive}
there is a close connection between deformations of Lie algebras and
the Lie algebra two-cohomology  $\H^2(\L,\L)$ with values in $\L$.
The following is obvious.
\begin{proposition}
Let $\L$ be an arbitrary Lie algebra over $\K$, then
$\L$ is infinitesimally rigid iff $\H^2(\L,\L)=0$.
\end{proposition}

\begin{theorem}\label{T:fd0}
(\cite{Ger}, \cite{NijRich1}, \cite{NijRich2}) 
Let $\L$ be a finite-dimensional Lie algebra.
If $\H^2(\L,\L)=0$ then the Lie algebra $\L$ is rigid in 
all the senses introduced above.
\end{theorem}
In particular, as infinitesimal rigidity is equivalent
to $\H^2(\L,\L)=0$, it already implies in finite dimension already rigidity in the
formal, geometric, (and analytic) sense.
This is definitely not true in infinite dimension, as we showed in
\cite{FiaSchlV} by exhibiting a deformation over the affine line 
$\C[t]$ of the Witt algebra which is non-trivial in every neighbourhood of
$t=0$, despite the fact that the Witt algebra is infinitesimally 
and formally rigid.
In \refS{elliptic} we will construct deformations of the 
current algebras $\gb$ 
associated to finite dimensional complex Lie algebras $\g$
which are neither geometrically nor analytically rigid despite
that for $\g$ simple the Lie algebra $\gb$ is 
formally rigid \cite{LeRoRid}.

Nevertheless, also in infinite dimensions we have  the following
theorem.
\begin{theorem}\label{T:formrig}
(\cite{Ger}, \cite{NijRich1}, \cite{NijRich2}) 
Let $\L$ be an arbitrary Lie algebra.
If  $\H^2(\L,\L)=0$
then $\L$ is formally rigid.
\end{theorem}
The theorem  is also a corollary of the result of 
 Fialowski \cite{Fiaproc} on the
existence of versal formal families which will be recalled in 
the next subsection.

\bigskip

For completeness 
and to contrast this with the infinite dimensional case 
let us quote the following results for
the {\bf finite-dimensional} case.
\begin{theorem}\label{T:thmrauch}
(\cite{Rau})
Let $\L$ be a finite-dimensional Lie algebra corresponding to the 
point $\mu\in Lalg_n$ then:
\newline
(a) The Zariski tangent space of the scheme $Lalg_n$ at the
point $\mu$ can be naturally identified with $\mathrm{Z}^2(L,L)$.
\newline
(b)  The Zariski tangent space  of the
$\mathrm{GL}(n)$ orbit of $\mu$ 
(considered as  reduced scheme) can be 
 naturally identified with $\mathrm{B}^2(L,L)$.
\end{theorem}
From \refT{thmrauch} and \refT{fd0}, 
and the fact  that $\H^2(\L,\L)= 0$
implies that $\mu$ is a non-singular point 
\cite{Ger}, \cite{NijRich2},
the following is a consequence.
\begin{theorem}\label{T:fdrig}
Let $\L$ be a  finite-dimensional Lie algebra.
Then  $\H^2(\L,\L)=0$ iff $\L$ is rigid and 
is a nonsingular Lie algebra.
Moreover, in this case all the above definitions of
rigidity coincide.
\end{theorem}
Richardson \cite{Rich} has given an example of a finite-dimensional
Lie algebra with $\H^2(\L,\L)\ne 0$, which is nevertheless 
rigid in the orbit sense. 
Clearly, the $\mu$ corresponding to the Lie structure of this $\L$
is a singular point in  $Lalg_n$.

\subsection{Universal and versal deformations}
$ $

As explained above and will be shown further down in this article, 
the strong relation between cohomology spaces and deformations
over non-local rings breaks down in infinite dimension.
It might be considered as a astonishing result that a  tight
connection still exists on the formal level.

\begin{proposition}\label{P:infuni}
(\cite{FiaFuMini})
Assume that  $\dim \H^2(\L,\L) < \infty$, then there exists a universal
infinitesimal deformation $\eta_\L$ of the Lie algebra $\L$ with
base $B= \K\oplus \H^2(\L,\L)^*$, where 
the second summand is the dual of $\H^2(\L,\L)$ equipped with 
the zero multiplication, i.e.
$$
(\alpha_1,h_1)\cdot(\alpha_2,h_2) = (\alpha_1\alpha_2,\alpha_1h_2 +
\alpha_2 h_1).
$$
\end{proposition}
This means that for every infinitesimal deformation $\lambda$ of the Lie
algebra $\L$ with finite dimensional base $A$, there exists a
unique homomorphism $\phi: \K\oplus \H^2(\L,\L)^*\to A$ such that
$\lambda$ is equivalent to the push-out $\phi_*\eta_\L$.
\medskip

Although in general it is impossible to construct a universal 
formal deformation, there is a so-called versal element.

\begin{definition} (\cite{Fiasbor,Fiaproc})
A formal deformation $\eta$ of $\L$ parameterized by a complete local
 algebra
$(B,m_B)$ is called {\it versal} if for every deformation $\lambda$,
parameterized by a complete local algebra $(A,m_A)$, there is a morphism
 $f: B\to A$ such that
\smallskip

\quad 1) The push-out $f_*\eta$ is equivalent to $\lambda$.
\smallskip

\quad 2) If $A$ satisfies ${m_A}^2=0$, then $f$ is unique.
\end{definition}

\begin{theorem}\label{T:formver}(\cite{Fiaproc},\cite[Thm. 4.6]{FiaFuMini})
Assume that  $\dim \H^2(\L,\L) < \infty$.
\newline
(a) There exists a  versal formal  deformation of $\L$.
\newline
(b) The base of the versal formal deformation is formally embedded into
$\H^2(\L,\L)$, i.e. it can be described in $\H^2(\L,\L)$ by
a finite system of formal equations.
\end{theorem}

Hence if $\H^2(\L,\L)=0$, every formal deformation will be equivalent
to the trivial one, see \refT{formrig}.
\subsection{Deformations  of commutative algebras and Harrison cohomology}
\label{S:harr}
$ $

In this article we will deform current algebras 
by deforming  associative and
commutative algebras
in a geometric way.
The corresponding cohomology theory of such deformations is the
Harrison cohomology \cite{Har}.

Let $\A$ be an associative and commutative algebra over $\K$.
We only need here the space $\H^2_{Harr}(\A,\A)$. Recall its
definition.
The two-cocycles are bilinear maps $F:\A\times \A\to\A$  such that
\begin{equation}
\begin{gathered}
F(a,b)=F(b,a),\qquad a,b\in\A,
\\
\delta_2 F(a,b,c):=aF(b,c)-F(ab,c)+F(a,bc)-F(a,b)c=0,\qquad a,b,c\in\A.
\end{gathered}
\end{equation}
A two-cycle is a coboundary if there exists a linear map
$\phi:\A\to\A$ such that 
\begin{equation}
F=\delta_1\phi(a,b)=
a\phi(b)-\phi(ab)+\phi(a)b,\qquad a,b\in\A.
\end{equation}
Note that  $\H^2_{Harr}(\A,\A)$ will be a subspace of the 
the Hochschild cohomology space 
 $\H^2_{Hoch}(\A,\A)$.

\section{Krichever-Novikov algebras}
\label{S:kn}
\subsection{The algebras with their almost-grading}
$ $

Algebras of Krichever-Novikov type are generalizations of the Virasoro
algebra, current algebras and all their related algebras.
Let $M$ be a compact Riemann surface of genus $g$, or in terms
of algebraic geometry, a smooth projective curve over $\C$.
Let $N,K\in\N$ with $N\ge 2$ and $1\le K<N$. Fix
$$
I=(P_1,\ldots,P_K),\quad\text{and}\quad
O=(Q_1,\ldots,Q_{N-K})
$$
disjoint  ordered tuples of  distinct points (``marked points'',
``punctures'') on the
curve. In particular, we assume $P_i\ne Q_j$ for every
pair $(i,j)$. The points in $I$ are
called the {\it in-points}, the points in $O$ the {\it out-points}.
Sometimes we consider $I$ and $O$ simply as sets and denote 
$A=I\cup O$ as a set.

Here we will need
the following algebras.
Let $\A$ be the associative algebra of those meromorphic functions 
on $M$ which are holomorphic outside the set of points $A$ with point-wise
multiplication.
Let $\L$ be the Lie algebra of meromorphic vector fields
which are holomorphic outside of $A$ with the usual Lie bracket of
vector fields.
The algebra $\L$ is called the {\it vector field algebra of 
Krichever-Novikov type}.
They were introduced and their structure was studied by 
Krichever and Novikov \cite{KNFac}.
The corresponding generalization to the multi-point case was done in
\cite{SchlDiss}, \cite{Schlkn}, \cite{Schleg}, \cite{Schlce}.
Obviously, both $\A$ and $\L$ are infinite dimensional algebras.

Furthermore we will need the {\it higher-genus, multi-point current
algebra of Krichever-Novikov type}.
We start with 
$\g$  a complex finite-dimensional Lie algebra and  endow
the tensor product $\Gb=\g\otimes_\C \A$ with the Lie bracket
\begin{equation}
[x\otimes f, y\otimes g]=[x,y]\otimes f\cdot g,
\qquad  x,y\in\g,\quad f,g\in\A.
\end{equation} 
The algebra  $\Gb$ is the higher genus current algebra.
It is an infinite dimensional Lie algebra and might be considered 
as the Lie algebra of $\g$-valued meromorphic functions on the
Riemann surface with poles only outside of $A$.

\medskip
The classical genus zero and $N=2$ point case is give by the 
geometric data
\begin{equation}\label{E:class}
M=\P^1(\C)=S^2, \quad I=\{z=0\},\quad 
 O=\{z=\infty\}.
\end{equation}
In this case the algebras are the well-known algebras of 
Conformal Field Theory (CFT).
For the function algebra
we obtain
$\A=\C[z^{-1},z]$, the algebra of Laurent polynomials.
The vector field algebra $\L$ is the Witt algebra
generated by $l_n=z^{n+1}\frac {d}{dz}$, $n\in \Z$ with
Lie bracket $[l_n,l_m]=(m-n)l_{n+m}$, and the current algebra $\Gb$ is
the standard current algebra
$\gb=\g\otimes \C[z^{-1},z]$  with
Lie bracket
\begin{equation}
[x\otimes z^n, y\otimes z^m]=[x,y]\otimes z^{n+m}
\qquad  x,y,,\quad n,m\in\Z.
\end{equation}

\medskip

For infinite dimensional algebras, modules
and their representation theory a graded structure is usually
 important to obtain structure results.
In fact, to define certain types of representations which are
fundamental,
 e.g. in CFT like the
highest weight representations, 
one uses grading.

In the classical situation the algebras are obviously graded by
taking as degree $\deg l_n:=n$ and $\deg x\otimes z^n:=n$.
For higher genus there is usually no grading.
But it was observed by Krichever and Novikov 
in the two-point case that a weaker
concept, an almost-graded structure, will be enough to develop an
interesting  theory of representations (Verma modules, etc.).
\begin{definition}\label{D:almgrad}
Let $\A$ be an (associative or Lie) algebra admitting a direct
decomposition as vector space $\ \A=\bigoplus_{n\in\Z} \A_n\ $.
The algebra $\A$ is called an {\it almost-graded}
algebra if (1) $\ \dim \A_n<\infty\ $ and (2)
there are constants $R$ and  $S$ such that 
\begin{equation}\label{E:eaga}
\A_n\cdot \A_m\quad \subseteq \bigoplus_{h=n+m+R}^{n+m+S} \A_h,
\qquad\forall n,m\in\Z\ .
\end{equation}
The elements of $\A_n$ are called {\it homogeneous  elements of degree $n$}.
\end{definition}
By exhibiting a special basis,
for the multi-point situation such an almost grading was introduced in
\cite{SchlDiss,Schlkn, Schleg,Schlce}.
Essentially, this is done by fixing the order of the basis elements
at the points in $I$ and
$O$ in a complementary way to make them unique.
In this way we obtain e.g. for the function algebra (resp. for the vector 
field algebra)
basis  elements
$A_{n,p}$ (resp. $e_{n,p}$) with $n\in\Z$ and $p=1,\ldots,K$.
As definition for the degree we take  ($x\in\g$)
\begin{equation}\label{E:grad}
\deg(e_{n,p}):=
\deg(A_{n,p}):=
\deg(x\otimes A_{n,p}):=n.
\end{equation}
In the following we will give an explicit description of the basis
elements
for those  genus zero and one situation we need.
Hence, we will not recall their general definition but only refer to
the above quoted articles.
\begin{proposition}
\cite{SchlDiss,Schlce}
With respect to the grading introduced by \refE{grad} the 
algebras $\L,\A$, and $\Gb$ 
are almost-graded.
The almost-grading depends on the splitting $A=I\cup O$.
\end{proposition}

\subsection{Central extensions}
\label{SS:cext}
$ $

In the construction of infinite dimensional representations of
these algebras with certain desired properties
(generated by a vacuum, irreducibility, unitarity, etc.)
one is typically forced to ``regularize'' a ``naive'' action
to make it well-defined. Important examples in CFT are  
the fermionic Fock space representations which are
constructed by taking semi-infinite forms of a fixed weight.

From the mathematical point of view,
with the
help of a prescribed  procedure one modifies the action to 
make it well-defined, but on the other hand, 
accepting that the modified action 
in compensation 
will be only  a projective Lie action. Such projective actions are
honest Lie actions  for a suitable centrally extended algebra.
In the classical case they are well-known.
The unique non-trivial (up to equivalence and rescaling) central extension
of the Witt algebra is the Virasoro algebra.
For the current algebra $\g\otimes \C[z^{-1},z]$ if $\g$ is a simple Lie
algebra 
with Cartan-Killing form $\beta$, it is the corresponding affine Lie algebra 
$\gh$ (or, untwisted affine Kac-Moody algebra):
\begin{equation}\label{E:classcur}
[x\otimes z^n, y\otimes z^m]=[x,y]\otimes z^{n+m}-\beta(x,y)\cdot 
n\cdot\delta_{m}^{-n}\cdot t,
\quad [t,\gh]=0,
\quad x,y\in\g,\ n,m\in\Z.
\end{equation}

For the extension to higher genus and many points the objects have
to be ``geometrized''.
First recall that for a Lie algebra $\V$ central extensions are 
classified (up to equivalence) by the second Lie algebra cohomology
$\H^2(\V,\C)$
of $\V$ with values in the trivial module $\C$.
A bilinear form $\psi:\V\times\V\to\C$ is called a Lie algebra 2-cocycle
iff $\psi$ is antisymmetric and fulfills the cocycle condition
\begin{equation}
0=d_2\psi(x,y,z):=
\psi([x,y],z)+
\psi([y,z],x)+
\psi([z,x],y).
\end{equation}
To obtain central extensions of $\Gb=\g\otimes \A$
we start with 
$\g$ being an arbitrary finite-dimensional Lie algebra
and $\beta$ a symmetric, invariant, bilinear form on it (not necessarily 
non-degenerate).
Invariance means  that we have 
$\beta([x,y],z)=\beta(x,[y,z])$ for all $x,y,z\in\g$.
We set 
$\Gh=\C\oplus\Gb$ as vector space and introduce 
the  bracket
\begin{equation}\label{E:daff}
[\widehat{x\otimes f},\widehat{y\otimes g}]=
\widehat{[x,y]\otimes (f g)}+\beta(x,y)\cins fdg \cdot t,\qquad
[\,t,\Gh]=0\ .
\end{equation}
Here  we
used for short $\widehat{x\otimes f}:=(0,x\otimes f)$,
 $t:=(1,0)$
and $C_S$ denotes  a cycle separating the points in $I$ from the points in
$O$.
\begin{proposition}[\cite{Schlaff}]
The term $\beta(x,y)\cins fdg$ is a Lie algebra 2-cocycle of $\Gb$ with values
in the trivial module $\C$. Hence,
the vector space
$\Gh$ with structure \refE{daff} is a 
Lie algebra which defines a central extension of $\Gb$.
\end{proposition}
These algebras are called {\it higher genus (multi-point)
affine algebras (of Krichever-Novikov type)}.
In the classical situation \refE{class} we obtain back \refE{classcur}.
In the following we will use the term {\it classical current algebra}
to denote the current algebra \refE{classcur}. In the new
terminology it is a genus zero and two-point current algebra.

Note that the cocycle can  be calculated as
\begin{equation}\label{E:cocycex}
\beta(x,y)\cdot \cins fdg=\beta(x,y)\cdot\sum_{k=1}^K\res_{P_k}(fdg)
=-\beta(x,y)\cdot\sum_{l=1}^{N-K}\res_{Q_l}(fdg).
\end{equation}

As in general a non-vanishing $\beta$ is not uniquely given 
(up to rescaling), the central extension $\Gh$ will depend
on it.
Even in the case when $\g$ is a simple Lie algebra, which implies
that there is essentially only one non-vanishing form $\beta$, the
Cartan-Killing form, for the higher genus or/and multi-point situation
$\H^2(\Gb,\C)$ will be more than one-dimensional
(e.g. take another path of integration  which is not
homologous to $C_S$).
But the following is shown in 
\cite{Schlaff}.
\begin{theorem}\label{T:localext}
Let $\g$ be a simple finite-dimensional Lie algebra, then up
to equivalence and 
multiplication with a scalar  there is a unique
non-trivial 2-cocycle class for 
the current algebra $\Gb=\g\otimes\A$ which has a ``local'' representative,
i.e. which allows to extend the almost grading 
of $\Gb$ to the central extension $\Gh$
defined by the local representative by assigning to the central element
$t$ a  degree.
This class is given by \refE{daff} with $\beta$ being the Cartan-Killing
form.
\end{theorem}
Also, in \cite{Schlaff} a thorough treatment for the 
 case when $\g$ is semi-simple or even reductive 
can be found.
Corresponding uniqueness results for almost-graded  central extensions
of $\A$ and $\L$ are shown in \cite{Schlloc}.

%
\section{Current algebras for the elliptic curve case}\label{S:elliptic}
In the following we will construct global deformations of the
classical current algebra $\gb=\g\otimes\C[z^{-1},z]$ 
and its central extension $\gh$.
First we construct 
a family of associative algebras which contain the algebra of 
Laurent polynomials 
$\C[z^{-1},z]$
as special element.
The deformation family for the current algebra will be obtained by 
tensoring $\g$ with this family.
These families will be of geometric origin.
More precisely, each non-special element in the family will be a
current
algebra of Krichever-Novikov type 
for the genus one (i.e. the elliptic) case.
The construction of these families
will not make any assumption about the finite-dimensional
Lie algebra.
Only later we will require $\g$ to be a simple
Lie algebra to contrast the existence of these global deformations
which are (geometrically) locally not equivalent to trivial
deformations, 
with the formal rigidity of $\gb$
\cite{LeRoRid}.

\subsection{The family of elliptic curves}
$ $

For the convenience of the reader
we will recall the geometric picture from
\cite{FiaSchlV}.
As we have geometric degenerations in mind, it is more convenient to pass 
from the complex analytic picture (i.e. the language of
Riemann surfaces) to the algebraic geometric picture 
(i.e. the language of curves).
Every compact Riemann surface of genus one corresponds to an
elliptic curve in the projective plane.
Recall that the elliptic curves can be given as sets of 
solutions  of the
polynomial equation
\begin{equation}\label{E:ellipp}
Y^2Z=4X^3-g_2XZ^2-g_3Z^3,\quad g_2,g_3\in\C,\quad
\text{with } \Delta:={g_2}^3-27{g_3}^2\ne 0.
\end{equation} 
Here $g_2$ and $g_3$ are parameterizing the individual curve and 
the condition $\Delta\ne 0$ assures that the curve will be
nonsingular.
Instead of \refE{ellipp} we can use the description
\begin{equation}\label{E:ellipe}
Y^2Z=4(X-e_1Z)(X-e_2Z)(X-e_3Z)
\end{equation}
with 
\begin{equation}\label{E:edef}
e_1+e_2+e_3=0,
\quad\text{and}\quad
\Delta=16(e_1-e_2)^2
(e_1-e_3)^2
(e_2-e_3)^2\ne 0.
\end{equation}
These presentations are related via
\begin{equation}
g_2=-4(e_1e_2+e_1e_3+e_2e_3),\quad
g_3=4(e_1e_2e_3).
\end{equation}
The elliptic modular parameter
classifying elliptic curves up to isomorphy is given as 
\begin{equation}
j=1728\;\frac {g_2^3}{\Delta}.
\end{equation}
We set 
\begin{equation}
B:=\{(e_1,e_2,e_3)\in \C^3\mid e_1+e_2+e_3=0,\quad
e_i\ne e_j\ \text{for}\ i\ne j\}.
\end{equation}
Inside the product $B\times\P^2$ we consider the family of elliptic
curves $\E$ over $B$ defined via \refE{ellipe},
with its obvious projection $\E\to B$.
The family can be extended to 
\begin{equation}
\Bh:=\{e_1,e_2,e_3)\in \C^3\mid e_1+e_2+e_3=0
\}.
\end{equation}
The fibers above $\Bh\setminus B$ are singular cubic curves.
Resolving the linear relation in $\Bh$ via
$e_3=-(e_1+e_2)$ we obtain a family over $\C^2$.

Consider the following complex lines in $\C^2$
\begin{equation}
D_s:=\{(e_1,e_2)\in\C^2\mid 
e_2=s\cdot e_1\},\quad  s\in \C,
\qquad
D_\infty:=\{(0,e_2)\in\C^2\}.
\end{equation} 
Set also
\begin{equation}
D_s^*=D_s\setminus\{(0,0)\} 
\end{equation}
for the punctured line.
Now 
\begin{equation}
B\cong\C^2\setminus (D_1\cup D_{-1/2}\cup D_{-2}).
\end{equation}
Note that above $D_1^*$ we have $e_1=e_2\ne e_3$,
 above $D_{-1/2}^*$ we have $e_2=e_3\ne e_1$,
and  above $D_{-2}^*$ we have $e_1=e_3\ne e_2$.
In all these cases we obtain a nodal cubic.
Every nodal cubic $E_N$ can be given as
\begin{equation}
Y^2Z=4(X-eZ)^2(X+2eZ)
\end{equation}
where $e$ denotes the value of the coinciding $e_i=e_j$  ($-2e$ is then 
necessarily the remaining one).
The singular point is the point $(e: 0:1)$. It is a node.

Above the unique common intersection point $(0,0)$ 
of all $D_s$ there is the
cuspidal cubic $E_C$:
\begin{equation} 
Y^2Z=4X^3.
\end{equation}
The singular point is $(0:0:1)$. 
In both cases the complex projective line is the desingularisation
of the singular curve.

In all cases (non-singular or singular) the point  
$\infty=(0:1:0)$ lies on the  
curves. It is the only intersection with the line at
infinity,  and is a non-singular point.
In  passing  to the affine plane in the following we will loose nothing.
The affine curve will be given as the solution set of
\begin{equation}\label{E:ellip}
Y^2=4(X-e_1)(X-e_2)(X-e_3).
\end{equation}
For the  curves above the points in $D_s^*$ we calculate $e_2=s e_1$ and 
$e_3=-(1+s)e_1$ (resp. $e_3=-e_2$ if $s=\infty$).
Due to the homogeneity,  the modular parameter $j$ 
for the  curves above $D_s^*$ 
will be constant along the line.
In particular, the curves in the family lying above $D_s^*$ will
be isomorphic.
Their modular parameter calculates to
\begin{equation}
j(s)=1728\;\frac {4(1+s+s^2)^3}{(1-s)^2(2+s)^2(1+2s)^2},
\quad j(\infty)=1728.
\end{equation}
\subsection{The family of current algebras}
$ $

First we define a family of function algebras $\A$ on 
these elliptic curves.
We introduce the points where poles are allowed.
For our purpose it is enough to consider two marked points.
More marked points are considered in \cite{SchlDeg}, \cite{RDS}.
We will always put one marking to $\infty=(0:1:0)$ and the 
other one to the point with the affine coordinate $(e_1,0)$.
These markings define two sections of the family $\E$ over $\Bh\cong
\C^2$.
With respect to the group structure on the elliptic curve
given by $\infty$ as the neutral element (the first marking), 
the second marking chooses a two-torsion point.
All other choices of two-torsion points will yield isomorphic 
situations.
\begin{proposition}\label{P:functalg}
For every elliptic curve $E_{(e_1,e_2)}$ 
over $(e_1,e_2)\in \C^2\setminus (\Dex)$ 
the associative  algebra $\A_{(e_1,e_2)}$ of functions on $E_{(e_1,e_2)}$ 
has a basis $\{A_n$, $n\in\Z\}$ such that   the algebra structure is
given as
\begin{equation}\label{E:structwo}
A_n\cdot A_m=
\begin{cases}
A_{n+m},& \text{for $n$ or $m$ even},
\\
A_{n+m}+3e_1 A_{n+m-2}
\\
\qquad +(e_1-e_2)(2e_1+e_2) A_{n+m-4},
& \text{for $n$ and $m$ both odd}.
\end{cases} 
\end{equation}
By setting $\deg(A_n):=n$, we obtain an almost-grading.
\end{proposition}
\begin{proof}
In \cite{SchlDeg} it was shown that a basis 
of the corresponding Krichever-Novikov function algebra  $\A$ is given by
\begin{equation}
A_{2k}:=(X-e_1)^{k},
\qquad
A_{2k+1}:=1/2 Y(X-e_1)^{k-1} \qquad k\in\Z.
\end{equation}
The calculation of the structure constants 
is straightforward. In the case when both $n$ and $m$ are odd, one
replaces $Y^2$ by 
$4(X-e_1)(X-e_2)(X-e_3)$
and uses  $e_1+e_2+e_3=0$.
The almost-grading is obvious.
\end{proof}
The algebras in \refP{functalg}  defined with the
structure \refE{structwo} make sense also for the points
$(e_1,e_2)\in D_1\cup D_{-1/2}\cup D_{-2}$.
Altogether this defines a two-dimensional family of  algebras
parameterized over the affine plane $\C^2$, or  described
algebraically, over the polynomial algebra $\C[e_1,e_2]$.
For $(e_1,e_2)=(0,0)$ we obtain the algebra of
Laurent polynomials, if we make the identification $A_n=z^n$.

By tensoring this family with $\g$ we get the following. 
\begin{theorem}\label{T:afffam}
Let $\g$ be a finite-dimensional Lie algebra, and
$V$ a vector space with basis $\{A_n\mid n\in\Z\}$ then
for given complex values $e_1,e_2$
the vector space $\g\otimes V$  carries
the structure of a Lie algebra given by
\begin{equation}\label{E:structaff}
[x\otimes A_n, y\otimes A_m]
=\begin{cases}
[x,y]\otimes A_{n+m},& \text{for $n$ or $m$ even},
\\
[x,y]\otimes A_{n+m}+3e_1 [x,y]\otimes A_{n+m-2}
\\
\qquad +(e_1-e_2)(2e_1+e_2) [x,y]\otimes A_{n+m-4},
& \text{for $n$ and $m$ both odd}.
\end{cases} 
\end{equation}
Here $x$ and $y$ are elements of $\g$.
In particular, these algebras define a two-parameter family
of deformations of the current algebra $\gb$, such that
$\gb$ corresponds to the point $(0,0)$ and the
algebra over $(e_1,e_2)\in B$ corresponds to the 
elliptic affine algebra $\Gb_{(e_1,e_2)}$ fixed by the geometric
data.
\end{theorem}
Further down we will identify the algebras over the other places.
But first we study 
the family of algebras  
obtained by taking  as base variety the (affine) line $D_s$ (for any $s$).
\begin{corollary}
For every $s\in\C\cup\{\infty\}$ the
families of \refT{afffam} define by restriction  one-dimensional
families of Lie algebras $\Gb^s_{(e)}$ which are deformations
of the classical current algebra $\gb$ (corresponding to $e=0$)
over the affine
line $\C$, resp. the algebra $\C[e]$. 
For  $s\ne\infty$ the family is given by
\begin{equation}\label{E:structs}
[x\otimes A_n, y\otimes A_m]
=\begin{cases}
[x,y]\otimes A_{n+m},& \text{for $n$ or $m$ even},
\\
[x,y]\otimes A_{n+m}+3e [x,y]\otimes A_{n+m-2}
\\
\qquad +e^2(1-s)(2+s)[x,y]\otimes A_{n+m-4},
& \text{for $n$ and $m$ both odd}.
\end{cases} 
\end{equation}
For $s=\infty$ the family is given by
\begin{equation}\label{E:struct0}
[x\otimes A_n, y\otimes A_m]
=\begin{cases}
[x,y]\otimes A_{n+m},& \text{for $n$ or $m$ even},
\\
[x,y]\otimes A_{n+m}-e[x,y]\otimes A_{n+m-4},
& \text{for $n$ and $m$ both odd}.
\end{cases} 
\end{equation}
In the case when $s\ne 1,-1/2$ or $2$ the 
algebras $\Gb^s_{(e)}$ are elliptic current algebras.
In any case for fixed $s$ we have  $\Gb^s_{(e)}\cong \Gb^s_{(e')}$
as long as both $e,e'\ne 0$.
\end{corollary}
\begin{proof}
If  $s\ne\infty$ then 
$e_2=s e_1$ and \refE{structaff}
reduces to \refE{structs} if we set $e:=e_1$.
For $s=\infty$ we have   $e_1=0$ and obtain
\refE{struct0} if we set $e:=e_2^2$.
If we rescale the basis elements  $A_n^*=(\sqrt{e})^{-n}A_n$ 
(for $s\ne\infty$),
 we obtain for $e \ne 0$
always the algebra with $e=1$ in our structure equations.
For $s=\infty$ a rescaling  $(\sqrt[4]{e})^{-n}A_n$
will do  the same (for $e\ne 0$).
Hence we see  that for fixed $s$ in all cases the algebras will be isomorphic 
above every point in $D_s$, as long as
we are not above $(0,0)$.
\end{proof}
The generic isomorphy class of the family $\Gb^s_{(e)}$ can be
given by taking $\Gb^s_{(1)}$. To avoid cumbersome notation we
often denote this isomorphy type just by  $\Gb^s$.
But we have to keep in mind that 
by \refP{curnon} the algebra $\Gb^s$ will not be isomorphic
to  the special element $\Gb^s_{(0)}=\gb$.

Clearly, if we do the same kind of restrictions to
$D_s$ for the families of associative algebras 
$\A_{(e_1,e_2)}$, we obtain similar 
one-dimensional algebraic families of commutative and associative algebras
which deform the algebra of Laurent polynomials.
We will denote these families by $\A^s_{(e)}$.
Again for fixed $s$ all algebras over $e\ne 0$ will be isomorphic and
we will denote this isomorphy type simply by $\A^s$.
\refP{assnon} will show that they are not isomorphic to the
special element, the algebra of Laurent polynomial.

\subsection{The three point case for $\P^1$}
$ $

There is another geometric family of algebras around. Its
geometric picture is the three point situation for the Riemann
sphere (the projective line).
Because we will need these algebras  later on anyhow, we will also give  
this family.
The geometric situation is $M=\P^1(\C)$, $I=\{\a,-\a\}$ and 
$O=\{\infty\}$, $\alpha\ne 0$.
Let us denote the corresponding Krichever-Novikov function algebra by
$\V_{(\alpha)}$.
\begin{proposition}
For every $\alpha\in\C^*$ the algebra $\V_{(\alpha)}$ has  a basis
$\{A_n\mid n\in\Z\}$ such that its structure is given by
\begin{equation}\label{E:strucptwo}
A_n\cdot A_m=
\begin{cases}
A_{n+m},& \text{for $n$ or $m$ even},
\\
A_{n+m}+\alpha^2A_{n+m-2},
& \text{for $n$ and $m$ both odd}.
\end{cases}
\end{equation}
By setting  $\deg(A_n):=n$ the algebra 
becomes an almost-graded algebra.
\end{proposition}
\begin{proof}
In \cite{SchlDeg} it was shown that a basis 
of  (the vector space) $\V_{(\alpha)}$ is given by
\begin{equation}
A_{2k}:=(z-\a)^k(z+\a)^k,
\qquad
A_{2k+1}:=z(z-\a)^{k}(z+\a)^{k}, \qquad k\in\Z.
\end{equation}
Here $z$ is the quasi-global coordinate on $\P^1(\C)$. 
The structure follows by direct calculations.
\end{proof}
Again, the algebra structure makes also sense for 
$\alpha=0$. In this case one obtains the algebra of Laurent polynomials.

For the corresponding current algebra 
$\mathcal{Y}_{(\alpha)}:=\g\otimes \V_{(\alpha)}$
we obtain the following.
\begin{theorem}
There exists a one-dimensional family of current algebras 
$\mathcal{Y}_{(\alpha)}$ with Lie structure
\begin{equation}\label{E:three}
[x\otimes A_n, y\otimes A_m]
=\begin{cases}
[x,y]\otimes A_{n+m},& \text{for $n$ or $m$ even},
\\
[x,y]\otimes A_{n+m}+\alpha^2[x,y]\otimes A_{n+m-2},
&\text{for $n$ and $m$ both odd},
\end{cases} 
\end{equation}
which is a deformation of the classical current algebra
$\gb$.
Moreover, $\mathcal{Y}_{(\alpha)}$  is equivalent to the families
$\Gb^{1}_{(e)}$ and $\Gb^{-2}_{(e)}$.
\end{theorem}
In \refS{degen} we will give a geometric explanation for the
latter identification.

\subsection{Isomorphy question}
$ $

In this subsection we will show that the current algebras 
$\Gb^{(s)}$ are not isomorphic to the classical current algebra
$\gb$ if $\g$ is semi-simple.
Hence, 
the families introduced are non-trivial algebraic-geometric deformations
for the  (formally) rigid  classical current algebra
$\gb$ (in case when $\g$ is simple).
\begin{proposition}\label{P:assnon}
The algebras $\A^{s}$  are not
isomorphic to the algebra of Laurent polynomials.
\end{proposition}
\begin{proof}
The involved algebras are the algebras of meromorphic functions
on  projective curves with at least one point removed.
Such curves are  affine curves and these algebras are
the affine coordinate algebras (i.e the algebras  of
regular  functions) of these curves. 
The isomorphy class of the coordinate algebra 
is uniquely given by the isomorphy class of the affine curve,
e.g. see \cite{Harts}.
Any algebra $\A^{s}$ for  $s\ne 1,-2,-1/2$ corresponds to an
elliptic curve with two points removed, 
both algebras $\A^{1}$ and $\A^{-2}$ correspond 
to the projective line with three points removed,
and $\A^{-1/2}$ corresponds to the nodal cubic with two (nonsingular)
points removed.
But the algebra of Laurent polynomials  
corresponds to the  projective line with two points removed.
The affine curves are not isomorphic
as can easily be seen from  the fact that their fundamental groups 
are different.
Hence, also the algebras are not isomorphic.
\end{proof}
\begin{proposition}\label{P:curnon}
Let $\g$ be a semi-simple finite-dimensional Lie algebra, 
and $\mathcal{A}$ and $\mathcal{B}$ 
two associative, commutative algebras (with units).
If the current algebras $\g\otimes \mathcal{A}$ and 
$\g\otimes \mathcal{B}$ are isomorphic as Lie algebras then
$\mathcal{A}$ and $\mathcal{B}$ are isomorphic as associative algebras.
\end{proposition}
\begin{proof}
Let $\mathfrak g$ be a semi-simple finite-dimensional Lie algebra, and
$P:\mathfrak g\otimes\mathcal{A}\to :\mathfrak g\otimes\mathcal{B}$ 
a Lie isomorphism of the current algebras.
The  Lie algebra $\mathfrak g$ 
admits a  $sl(2)$ subalgebra 
and hence also $\mathfrak g\otimes\mathcal{A}$ 
admits a  $sl(2)\otimes \mathcal{A}$ subalgebra. 
By restriction we obtain a Lie isomorphism
$P:sl(2)\otimes\mathcal{A} \to P(sl(2)\otimes\mathcal{A})$.
Via $sl(2)\cong sl(2)\otimes 1_A$ ($1_A$ the unit in $\mathcal{A}$),
$sl(2)$ is a subalgebra of $sl(2)\otimes\mathcal{A}$.
Denote by $h,e,f$ the standard generators of $sl(2)$, i.e.
\begin{equation}
h=
\begin{pmatrix}
1&0\\
0&-1
\end{pmatrix},
\ 
e=
\begin{pmatrix}
0&1\\
0&0
\end{pmatrix},
f=
\begin{pmatrix}
0&0\\
1&0
\end{pmatrix},
\quad
[h,e]=2e,\  [h,f]=-2f,\  [e,f]=h.
\end{equation}
The image $P(sl(2)\otimes 1_A)$ is isomorphic to $sl(2)$.
In particular, $P(h\otimes 1_A)$ will be mapped to a basis of its
Cartan subalgebra. After applying an inner automorphism,
we can assume that   $P(h\otimes 1_A)=h\otimes a_1$ with
$a_1\in B$, $a_1\ne 0$.
Let 
\begin{equation}
P(e\otimes 1_A)=h\otimes b_1+e\otimes b_2+f\otimes b_3,
\quad
P(f\otimes 1_A)=h\otimes c_1+e\otimes c_2+f\otimes c_3,
\end{equation}
with $b_i,c_i\in\mathcal{B},i=1,2,3$.
Using the structure equations above we see that there exist only two
solutions:
(A) $a_1=1$, $b_2=\alpha$ with an invertible element 
$\alpha\in \mathcal{B}$, $c_3=\alpha^{-1}$, and 
all other elements equal zero; and
(B)  $a_1=-1$, $b_3=\alpha$ with an invertible element 
$\alpha\in \mathcal{B}$,  $c_2=\alpha^{-1}$, 
and all other elements equal zero.
After an inner automorphism (given by
$\begin{pmatrix}
0&1\\ -1&0 \end{pmatrix}$) the solution (B) is transferred to
solution (A).
Hence, we can assume that 
$P$ has the property 
\begin{equation}
P(h\otimes 1_A)=h\otimes 1_B,\quad
P(e\otimes 1_A)=e\otimes \alpha,\quad
P(f\otimes 1_A)=f\otimes \alpha^{-1},\quad 
\text{$\alpha\in \mathcal{B}$ invertible}.
\end{equation}
We can decompose $P(sl(2)\otimes \mathcal{A})$ into the weight spaces
under the action of $h\otimes 1_B=P(h\otimes 1_A)$.
The weight spaces 
$h\otimes \mathcal{B}$, $e\otimes \mathcal{B}$,
and  $f\otimes \mathcal{B}$ correspond to the weights $0,2$, and
$-2$ respectively.
As
\begin{equation}
[h\otimes 1_B,P(x\otimes a)]
=
[P(h\otimes 1_A),P(x\otimes a)]=
P([h\otimes 1_A,x\otimes a])=
P([h,x]\otimes a),
\end{equation}
for $a\in\mathcal{A}$,
we obtain for $x\in \{h,e,f\}$
\begin{equation}
P(h\otimes a)=h\otimes Q^h(a),\quad
P(e\otimes a)=e\otimes  Q^e(a),\quad
P(f\otimes a)=f\otimes  Q^f(a),
\end{equation}
with linear maps $Q^h,Q^e,Q^f:\mathcal{A}\to \mathcal{B}$.
We will show that the map $Q^h$ is an algebra isomorphism.
Let $a_1,a_2\in \mathcal{A}$. 
{F}rom $[e\otimes a_1,f\otimes a_2]=[e,f]\otimes (a_1a_2)=
h \otimes (a_1a_2)$ one obtains after applying $P$ that
$Q^h(a_1a_2)=Q^e(a_1)Q^f(a_2)$.
Next we consider $[h\otimes a,e\otimes 1_A]=[h,e]\otimes a
=2e\otimes a$. As $\alpha=Q^e(1)$, we obtain
after applying $P$ that 
$Q^e(a)=Q^h(a)\cdot \alpha$.
In the same way we get $Q^f(a)=Q^h(a)\cdot \alpha^{-1}$.
Hence, $Q^h(a_1\cdot a_2)=Q^h(a_1)\cdot Q^h(a_2)$ and $Q^h(1)=1$.
This implies that $Q^h$ is indeed an isomorphism 
of associative algebras $\mathcal{A}\to \mathcal{B}$. 
\end{proof}
The proof shows that \refP{curnon} is also true for the case that
$\mathfrak g$ is a reductive but non-abelian Lie algebra.
In fact, it is only required that $\mathfrak g$ admits a Lie subalgebra
isomorphic to $sl(2)$.

\begin{theorem}\label{T:families}
Let $\g$ be a finite-dimensional simple Lie algebra. Then the 
classical current algebra $\gb=\g\otimes\C[z^{-1},z]$ is formally
rigid, but neither geometrically nor analytically rigid.
Examples of nontrivial geometric  deformations over the affine line $\C$, 
resp. the algebra
$\C[e]$, are given by the current algebras 
$\Gb^s_{(e)}$ of 
Krichever-Novikov type  \refE{structs}, \refE{struct0} for
every $s\in\C\cup \{\infty\}$.
\end{theorem}
\begin{proof}
First note that as $\g$ is simple, the classical current algebra $\gb$ is
rigid by a result of Lecomte and Roger \cite{LeRoRid}.
But obviously for every $s$,  $\Gb^s_{(e)}$ defines a geometric
deformation of $\gb$ over $\C[e]$. In particular we get
$\Gb^s_{(0)}=\gb$.
As  by \refP{assnon} the function algebras are not
isomorphic, \refP{curnon} implies that the algebra 
$\Gb^s_{(e)}$ is also not isomorphic to
$\gb$ as long $e\ne 0$.
Hence, restricted to every (algebraic-geometrically or
even analytically) open neighbourhood of $e=0$, the family will not be 
equivalent to the trivial family.
By Definition \ref{D:rigidalg} and \ref{D:rigidana}, $\gb$ is neither
geometrically, nor analytically rigid.
\end{proof}

\subsection{Families of the centrally extended algebras}
\label{S:famext}
$ $

The above family of current algebras $\Gb_{(e_1,e_2)}$  
resp. the one-dimensional families 
$\Gb^s_{(e)}$,  
can be 
centrally extended to a family of algebras   $\Gh_{(e_1,e_2))}$, 
resp. $\Gh^s_{(e)}$, i. e. to families of 
higher genus affine Lie algebras.
This can be achieved by using the defining equation \refE{daff}.
\begin{theorem}\label{T:famcent}
Assume we have a finite-dimensional Lie algebra $\gb$ and a
symmetric invariant bilinear form $\beta$ on $\g$.
For the family  $\Gb_{(e_1,e_2)}$ of 
current algebras 
\refE{structaff} containing the classical current algebra
$\gb$ as special element a family of almost-graded central extensions
 $\Gh_{(e_1,e_2)}$, i.e. a family of higher genus
affine Lie algebras  taining the classical affine
Lie algebra $\gh$ (with respect to the form $\beta$)
as special element, is given via the geometric two-cocycle
\refE{daff}. It calculates as
\begin{equation}
\label{E:famcent}
\gamma(x\otimes A_n,y\otimes A_m)=
p(e_1,e_2)\cdot \beta(x,y)\cdot\cins A_ndA_m
\end{equation}
with 
\begin{equation}
\label{E:famcent1}
\cins A_ndA_m=
\begin{cases}
\qquad -n\delta_m^{-n},&\text{$n$, $m$ even},
\\
\qquad\qquad 0,&\text{$n$, $m$ different parity},
\\
-n\delta_m^{-n}+
3e_1(-n+1)\delta_m^{-n+2}+
\\
+(e_1-e_2)(2e_1+e_2)(-n+2)\delta_m^{-n+4}
,&\text{$n$, $m$ odd}.
\end{cases}
\end{equation}
Here $p(e_1,e_2)$ is an arbitrary polynomial in the variables
$e_1$ and $e_2$.
\end{theorem}
The proof involves residue calculus (see Equation \refE{cocycex})
and will be postponed to an appendix. 
As the cocycle values \refE{famcent1} vanishes if $0\le n+m\le 4$,
the centrally extended algebras are almost-graded
by setting $\deg t:=1$ (or any other fixed value).

Clearly, 
for $e_1=e_2=0$ we obtain the classical affine algebra $\gh$.
By restricting this two-dimensional family to the
lines $D_s$ we get one-dimensional families.
For $s\ne\infty$ this amounts to replacing
in the last term $(e_1-e_2)(2e_1+e_2)$ by
$e_1^2(1-s)(2+s)$.
In particular, over $D_1$ and $D_{-2}$ it will 
vanish. 

By the uniqueness result for almost-graded central extensions
(\refT{localext}) in the case that $\g$ is simple we obtain
\begin{corollary}\label{C:famext}
If $\g$ is a finite-dimensional simple Lie algebra
then $\beta$ is necessarily the Cartan-Killing form, and
$\Gh_{(e_1,e_2)}$ \refE{famcent} 
is up to equivalence the unique almost graded central extension
of the family of current algebras $\Gb_{(e_1,e_2)}$.
\end{corollary}
Based on  the classification results of \cite{Schlaff} in the
semi-simple
case,  \refC{famext} remains true if one replaces
the factor $p(e_1,e_2)\beta(x,y)$ by
$\sum_{i=1}^K p_i(e_1,e_2)\beta_i(x,y)$, where again $p_i$ are
polynomials, and $\beta_i$ are the Cartan-Killing forms 
of the $K$ simple summands of $\gb$.
A similar statement is true for the reductive case, now 
additionally requiring the defining cocycle to be $\L$-invariant
(see \cite{Schlaff} for the definition).

\section{Geometric degenerations}
\label{S:degen}
It might be quite instructive to identify geometrically
all the algebras corresponding
to the singular cubic situations.
Above,
by comparing the structure constants
 we identified $\Gb_{(0,0)}$ and
$\Gb^{(s)}$ for $s=1$ or $-2$ with algebras which
are related to the genus zero situation .
The  geometric scheme behind this is that in each case the
desingularisation
(or normalization) of the singular cubic is the projective line.
By pulling back functions on  the singular cubics we obtain 
functions on the desingularisation.
But not necessarily all functions on the desingularisation
will be obtained as pullbacks.
 One also has to keep track of the poles.
In \cite{SchlDeg} the situation is analysed in complete details.

Three different situations appear.
\begin{itemize}
\item[(I)]
All three $e_1,e_2$ and $e_3$ coincide. 
The normalization $e_1+e_2+e_3=0$  implies
necessarily that they are zero. We obtain the cuspidal cubic
where the singular point is also a point of possible poles.
Over the singular point there is only one point.
Hence we obtain an identification of the algebra 
$\A_{(0,0)}$ with the full algebra $\C[z,z^{-1}]$,
and furthermore the identification of the current algebra
$\Gb_{(0,0)}$ 
with the classical current algebra $\gb$.
\item[(II)]
If two of the $e_i$ coincide (but not all three), we obtain the
nodal cubic. Over the singular point we have two points on
the desingularisation. Let the desingularisation
be chosen such that the points $\alpha$ and $-\alpha$ lie above
the singular point. 
We have to distinguish two sub-cases.
\begin{itemize}
\item[(IIa)]
 Either $e_1=e_2\ne e_3$ or $e_1=e_3\ne e_2$, then the singular point
(the node) will become a possible pole.
This situation occurs if we approach points from $ D^*_{1}\cup D^*_{-2}$.
The algebra generated by the pullbacks will be the full function
algebra of the 3-point algebra $\V_{(\alpha)}$. 
Hence, we see geometrically the already remarked isomorphy.
On the level of the current algebra we get the identification of
$\mathcal{Y}_{(\alpha)}$ with $\Gb_{(e)}^{1}$ and 
$\Gb_{(e)}^{-2}$.
\item[(IIb)]
If $e_1\ne e_2=e_3$, then the point of a possible pole will remain
non-singular.
This appears if we approach a point of $D^*_{-1/2}$.
For the pullbacks of the functions it is now
necessary that they have the same value at the points
$\alpha$ and $-\alpha$.
Hence, all elements of the algebra generated by the pullbacks
will have the same property.
We describe this algebra in the following.  
\end{itemize}
\end{itemize}
\begin{proposition}\label{P:under}
The set of elements 
\begin{equation}
A_n:=
\begin{cases}
z^n,& \text{$n$ even},
\\
z^n-\alpha^2z^{n-2}=z^{n-2}(z^2-\alpha^2),
&\text{$n$ odd},
\end{cases}
\end{equation}
for $n\in\Z$ form a basis of the subalgebra $\W_{(\alpha)}$ of 
meromorphic functions on
$\P^1$ which are holomorphic  outside $0$ and $\infty$ and have
the same value at $\alpha$ and $-\alpha$.
The algebra structure is given by
\begin{equation}
A_n\cdot A_m=
\begin{cases}
A_{n+m},& \text{for $n$ or $m$ even},
\\
A_{n+m}-2\alpha^2 A_{n+m-2}
\\
\qquad +\alpha^4 A_{n+m-4},
& \text{for $n$ and $m$ both odd}.
\end{cases} 
\end{equation}
For this subalgebra we have
$\mathcal{W}_{(\alpha)}\cong \A^{-1/2}$ and it is not isomorphic to
$\C[z^{-1},z]$.
\end{proposition}
\begin{proof}
Obviously these elements lie in this subalgebra and form a basis of
the  subalgebra of sums of all even functions (without any
restriction)
and odd functions
which vanish at $\pm\alpha$.
Let $f$ be a function fulfilling the conditions
for being a member of $\W_{(\alpha)}$.
Decompose it into its symmetric and antisymmetric part
\begin{equation}
f=f_1+f_2, \quad
f_1(x)=1/2(f(x)+f(-x)), \quad
f_2(x)=1/2(f(x)-f(-x)).
\end{equation}
Obviously $f_1$ also fulfills the conditions, hence
$f_2$ too and we get 
$f_2(-\alpha)=f_2(\alpha)$.
Being an antisymmetric function this implies that 
$f_2$ has to vanish at  $\pm\alpha$, and
$f$ is a linear combination of the elements $A_n$.
By setting $\alpha=\mathrm{i}\sqrt{\frac{3e_1}{2}}$ we immediately see that
$\mathcal{W}_{(\alpha)}\cong \A^{-1/2}$.
Hence by \refP{assnon}, $\W_{(\alpha)}$ is not isomorphic to 
$\C[z^{-1},z]$.
\end{proof}
Clearly again the algebras $\mathcal{W}_{(\alpha)}$ are isomorphic for different
$\alpha\ne 0$.

The above mentioned identifications extend to the current algebras
and we obtain another one-dimensional 
algebraic-geometric deformation family $\mathcal{Z}_{(\alpha)}$ of the current algebra.
\begin{proposition}
In the two-parameter family \refE{structaff},  the current algebras
for the singular cubic cases 
are isomorphic as follows:
\begin{equation}
\begin{aligned}
\Gb_{(0,0)}&\cong \gb=\g\otimes \C[z,z^{-1}],
\\
\Gb^{1}\cong \Gb^{-2}&\cong \g\otimes \mathcal{V},
\\
\Gb^{-1/2}&\cong \g\otimes \mathcal{Z}.
\end{aligned}
\end{equation}
\end{proposition}
Recall that here we denoted the  family with the
same letter as the isomorphy type of the  generic member.

These considerations extend to  the centrally extended 
families of algebras introduced in \refS{famext}, see 
\refP{singext}.
\section{Cohomology classes of the deformations}
\label{S:cohomo}
Let $\Gb_t$ be a one-parameter deformation of the current algebra $\gb$
with Lie structure 
\begin{equation}\label{E:family}
[\alpha,\beta]_t=[\alpha,\beta]+
t^k\omega_0(\alpha,\beta)+
t^{k+1}\omega_{1}(\alpha,\beta)+\cdots,
\end{equation}
such that $\omega=\omega_0$ is non-vanishing.
As explained in \refS{intuitive}, the bilinear map $\omega$ will be an element
of $\mathrm{C}^2(\gb,\gb)$.

For the global families $\Gb^s_{(e)}$ over the 
affine line $D_s$ appearing in \refS{elliptic}
we obtain as first nontrivial contribution the following two-cocycles 
($s\ne\infty$):
\begin{equation}\label{E:cocycl1}
\omega(x\otimes A_n, y\otimes A_m)
=\begin{cases}
0,& \text{for $n$ or $m$ even},
\\
[x,y]\otimes A_{n+m-2}
& \text{for $n$ and $m$ both odd}.
\end{cases} 
\end{equation}
For $D_\infty$ we get
\begin{equation}\label{E:cocycl2}
\omega(x\otimes A_n, y\otimes A_m)
=\begin{cases}
0,& \text{for $n$ or $m$ even},
\\
[x,y]\otimes A_{n+m-4},
& \text{for $n$ and $m$ both odd}.
\end{cases} 
\end{equation}
Here we used   the  symbols $x\otimes A_n$ for the basis elements
$x\otimes z^n$ also in the classical case. 
It should be pointed out that the families
corresponding to the singular cubic situations are
included for $s=1,-2$ or $s=-1/2$.
\begin{proposition}\label{P:trivial}
Let $\g$ be an arbitrary finite dimensional Lie algebra.
The two-cocycles \refE{cocycl1} and \refE{cocycl2} for the algebra
$\gb$ 
are cohomologically trivial. Hence, considered as infinitesimal
deformations of $\gb$, the families $\Gb^s_{(e)}$ 
are trivial.
\end{proposition}
\begin{proof}
We verify directly that 
$\omega=d_1\eta$ for $\eta:\gb\to\gb$ defined by 
\begin{equation}
\eta(x\otimes A_n)=
\begin{cases}
0,&n\text{\ even},
\\
-\frac 12x\otimes A_{n-2},&n\text{\ odd,}
\end{cases}
\quad\text{resp.}
\quad
\begin{cases}
0,&n\text{\ even},
\\
-\frac 12x\otimes A_{n-4},&n\text{\ odd.}
\end{cases}
\end{equation}
Clearly, if $\g$ is a simple Lie algebra the vanishing
of the cohomology class follows from the formal rigidity of
$\gb$ (\cite{LeRoRid}).
\end{proof}
\begin{remark}
As an intermediate step we constructed 
(inside the category of commutative and associative algebras) nontrivial 
algebraic-geometric deformations
$\A^s_{(e)}$ of the 
algebra of Laurent polynomials $\C[z^{-1},z]$.
One sees here the same effect as in the Lie case.
The relevant cohomology theory is the Harrisson cohomology
(see \refS{harr}).
As   $\C[z^{-1},z]$ is the algebra of regular functions on 
the smooth affine curve $\P^1\setminus\{0,\infty\}$, it is 
a smooth affine algebra and its Harrisson two-cohomology vanishes
\cite[Thm.22]{Har}.
Again we get that,  despite the fact that this algebra is
infinitesimally
and formally rigid, there exist nontrivial local  geometric families.
See Kontsevich \cite{Kon} for his concept of semi-formal
deformations (related to filtrations
of certain type) to overcome this discrepancy.
Indeed, by the almost-gradedness of the involved families 
considered in this article, the families  $\A^s_{(e)}$ are semi-formal
deformations in his sense.
Further considerations in these directions have  to be postponed to
another article.

As the Harrison two-cohomology vanishes, the defining cocycle for
these families will be a coboundary. One verifies immediately that
it is the coboundary of the linear form
\begin{equation}
\eta( A_n)=
\begin{cases}
0,&n\text{\ even},
\\
\frac 12A_{n-2},&n\text{\ odd,}
\end{cases}
\qquad\text{resp.}
\quad
\begin{cases}
0,&n\text{\ even},
\\
\frac 12A_{n-4},&n\text{\ odd.}
\end{cases}
\end{equation}

\end{remark} 

\appendix
\section{Proof of \refT{famcent}}\label{S:proof}
In this appendix we will show that the geometrically defined
central extension \refE{structaff} of the family $\Gb_{(e_1,e_2)}$ 
has the form \refE{famcent} if expressed 
for  pairs of generators $x\otimes A_n$, $x\in\g$ and $A_n$,
$n\in\Z$ the (homogeneous) basis elements of $\A_{(e_1,e_2)}$.
It is enough
to show  the expression \refE{famcent1}
for $\cins A_ndA_m$ for 
every pair of basis elements 
of the function algebra $\A_{(e_1,e_2)}$.
As a remark aside, 
this  integral  defines central extensions
of the algebra of functions considered as abelian Lie algebras,
see \cite{Schlloc}.

As the integration is over a separating cycle $C_S$, the integral can
be calculated by calculating residues at either one of
 the possible poles,
$(e_1,0)$ or $\infty$. From the point $\infty$ the residue has to be
taken with a minus sign.
One possible way to calculate the residue is to change to
the complex-analytic picture.
This means we use the fact that $X$ corresponds to the
elliptic Weierstra\ss\ $\wp$-function  and $Y$ to its derivative $\wp'$.
They are doubly-periodic functions on the complex plane 
with respect to the lattice $\Gamma=\Z\oplus\tau \Z$
with $\tau\in\C$, $\mathrm{im\,}\tau>0$.
We use the variable $z$ for the complex variable in the plane.
With this identification the complex one-dimensional torus
$\C/\Gamma$ is analytically isomorphic to the 
(projective) elliptic curve.
The variable $\tau$ fixes the isomorphy class of the 
torus, the parameters $e_1,e_2,e_3$ the isomorphy class of
the elliptic curve.
Clearly, they are related. One relation which is important for
us is
\begin{equation}
\wp(\frac {1}{2})=e_1,\qquad
\wp(\frac {\tau}{2})=e_2,\qquad
\wp(\frac {\tau+1}{2})=e_3.
\end{equation}
Note also that the 
functions $\wp$ and $\wp'$ depend on $\tau$ and hence also 
on the parameters $e_1,e_2$ and $e_3$.

We recall the following well-know facts from the
theory of elliptic functions. 
The function $\wp$ fulfills the differential equation
\begin{equation}\label{E:diffequ}
(\wp')^2=4(\wp-e_1)(\wp-e_2)(\wp-e_3)=
4\wp^3-g_2\wp-g_3.
\end{equation}

We have
\begin{equation}\label{E:umrech}
e_1+e_2+e_3=0,\qquad
g_2=-4(e_1e_2+e_1e_3+e_2e_3),
\qquad
g_3=4(e_1e_2e_3).
\end{equation}
The function $\wp$ is an even meromorphic function 
with poles of order two at the
points of the lattice and holomorphic elsewhere.
The function $\wp'$ is an odd meromorphic function with poles of order 
three at the
points of the lattice  and holomorphic elsewhere.
It has zeros of order one at the points $1/2,\tau/2$ and $(1+\tau)/2$
and all its translates under the lattice.
The  zeros are of order one.

For the Laurent series expansion at $z=0$ we obtain
\begin{equation}\label{E:pexp}
\wp(z)=\frac 1{z^2}(1+\frac {g_2}{20}z^4+O(z^6)),
\qquad
\wp(z)-e_1=\frac 1{z^2}(1+e_1z^2+\frac {g_2}{20}z^4+O(z^6)).
\end{equation}
  
In terms of these functions the  basis elements $A_n$ can be expressed as
\begin{equation}
A_{2k}=(\wp-e_1)^{k},
\qquad
A_{2k+1}=\frac 12 \wp'\cdot (\wp-e_1)^{k-1} \qquad k\in\Z.
\end{equation}
Note that $(\wp-e_1)$ has  a pole of order two at $z=0$ and
a zero of order two at $z=1/2$.
It has a Laurent expansion in even powers of $z$. Hence the same
is true for $A_{2k}$.
Furthermore, 
$A_{2k+1}(z)=\frac 1{2k}\frac {d}{dz}A_{2k}(z)$.
This  allows us to determine the Laurent expansion of
$A_{2k+1}$ easily.
Its Laurent expansion  consists
of only odd powers of $z$.

Our first conclusion is  that
the differential $A_{2k+1}A'_{2l}dz$ has only  
even terms in its Laurent expansion. Hence, as the residue
is the coefficient of $\frac 1z$, there is no
residue, and the cocycle values 
evaluated for pairs of basis elements of different 
parity is  zero, as claimed \refE{famcent1}.

Due to the residue theorem (on the torus) the
total residue  of any differential has to vanish. Hence, there
will be only a residue  at one point if the differential has poles at both
points, $z=0$ and $z=1/2$.
In the even case, 
$A_{2k}A'_{2l}dz$,
the pole order at $0$ is $2(k+l)+1$, and at
$1/2$ it is $(-2k)+(-2l+1)=-2(k+l)+1$.
Hence, a  possibly non-vanishing value could appear
only if $l=-k$ (or, equivalently  $m=-n$).
In the odd case,  
$A_{2k+1}A'_{2l+1}dz$,
we have poles at $0$ of order 
$(2k+1)+(2l+2)=2(k+l)+3$, and  at $1/2$ of order
$(-2k+1)+(-2l+2)=-2(k+l)+3$.
Hence, there could only be a contribution if
$l=-k-1,-k$, or $-k+1$ corresponding to $m=-n,-n+2$, or $-n+4$.

To calculate the residue we have to consider the
Laurent expansion of the generators
and their derivatives.
Starting from \refE{pexp} we get
\begin{equation}
\begin{aligned}
A_{2k}
=& \frac 1{z^{2k}}
\left(1+ke_1z^2+k\big(\frac {k-1}{2}e_1^2+\frac {g_2}{20}\big)z^4+O(z^6),
\right)
\\
A_{2k}'=&\frac 1{z^{2k+1}}(-2k+O(z^2)),
\\
A_{2k+1}=\frac 1{2k}A_{2k}'=
&\frac 1{z^{2k+1}}
\left(-1+(-k+1)e_1z^2+(-k+2)\big(\frac {k-1}{2}e_1^2+\frac {g_2}{20}\big)
z^4+O(z^6),\right)
\\
A'_{2k+1}=&\frac 1{z^{2k+2}}
\bigg(2k+1+(-k+1)(-2k+1)e_1z^2+
\\
&\qquad\qquad
+(-k+2)(-2k+3)
\big(\frac {k-1}{2}e_1^2+\frac {g_2}{20}\big)
z^4+O(z^6),\bigg)
\end{aligned}
\end{equation}
For the even pairing, i.e. $n=2k$, $m=2l$,
we obtain in case $l=-k$
\begin{equation}
A_{2k}A'_{-2k}dz=\frac 1z(2k+O(z^2)).
\end{equation}
Hence the residue at $z=0$ is $2k=m$, and as it has to be taken 
negatively, we obtain
$\cins A_{n}dA_{-n}=-n$, as claimed in \refE{famcent1}.

For the odd pairing, i.e. $n=2k+1$, $m=2l+1$,
  we have to multiply the Laurent expansions of
$A_{2k+1}$ and $A'_{2l+1}$. The resulting Laurent series is
\begin{equation}
\frac 1{z^{2(k+l)+3}}
\left(a(k,l)+b(k,l)z^2+c(k,l)z^4+O(z^6)\right)dz,
\end{equation}
with $a(k,l), b(k,l)$ and $c(k,l)$ obtained by collecting the 
corresponding contributions to these  orders.
In particular, we obtain
\begin{equation}
a(k,l)=-(2l+1),\qquad
b(k,l)=e_1((-k+1)(2l+1)+l-1)
\end{equation}
Now the residue at $z=0$ can be easily calculated.
\newline
For $l=-k-1$ we obtain as residue 
$a(k,-k-1)=(2k+1)=n$.
\newline
For $l=-k$ the residue is
$b(k,-k)=-3e_1(-n+2)$.
\newline
For $l=-k+1$ the residue is
\begin{equation}
c(k,-k+1)=
-(e_1-e_2)(2e_1+e_2)(-n+4).
\end{equation}
For the last calculation we have to express $g_2$ by the $e_i$'s 
as given in \refE{umrech} and replace $e_3=-(e_1+e_2)$.
Finally, it should not be forgotten to change the sign to
obtain the cocycle values claimed in \refE{famcent}.

\medskip
Strictly speaking, this derivation via the elliptic functions
is valid above the points outside of the collection of lines  
$D_1\cup D_{-2}\cup D_{-1/2}$. But  being a cocycle is a closed 
condition, hence the bilinear form will also be a cocycle for the 
algebras not corresponding to elliptic algebras.
It will be obtained by specializing the values of $e_1$ and $e_2$.

Nevertheless, it is an easy exercise to calculate the cocycle
for the genus zero algebras using again residues.
In this case the functions are rational functions, hence the 
calculations will be easier.
For a convenient  reference, let us note 
here the form of the cocycle for the singular situation.
\begin{proposition}\label{P:singext}
(a) For the current algebra $\Gb_{(0,0)}=\gb$ the centrally
extended algebra is given by the cocycle 
($p\in\C$)
\begin{equation}
\gamma(x\otimes A_n,y\otimes A_m)=
p\cdot \beta(x,y)\cdot (-n)\delta_m^{-n}.
\end{equation}
(b) For the three-point genus zero family of current algebras 
$\mathcal{Y}_{(\alpha)}$
\refE{three}
the family of centrally extended algebras is given by
the cocycle
\begin{equation}
\gamma(x\otimes A_n,y\otimes A_m)=
p(\alpha)\cdot \beta(x,y)\cdot\cins A_ndA_m,
\end{equation}
with a polynomial $p$ in $\alpha$ and 
\begin{equation}
\cins A_ndA_m=
\begin{cases}
\qquad -n\delta_m^{-n},&\text{$n$, $m$ even},
\\
\qquad\qquad 0,&\text{$n$, $m$ different parity},
\\
-n\delta_m^{-n}+
\alpha^2(-n+1)\delta_m^{-n+2},&\text{$n$, $m$ odd.}
\end{cases}
\end{equation}
\newline
(c) For the family of current algebras $\mathcal{Z}_{(\alpha)}$  
defined in \refS{degen} 
the family of centrally extended algebras is given by the cocycle
\begin{equation}
\gamma(x\otimes A_n,y\otimes A_m)=
p(\alpha)\cdot \beta(x,y)\cdot\cins A_ndA_m,
\end{equation}
with   a polynomial $p$ in $\alpha$ and 
\begin{equation}
\cins A_ndA_m=
\begin{cases}
\qquad -n\delta_m^{-n},&\text{$n$, $m$ even},
\\
\qquad\qquad 0,&\text{$n$, $m$ different parity},
\\
-n\delta_m^{-n}-2\alpha^2
(-n+1)\delta_m^{-n+2}+
\\
\quad +\alpha^4(-n+2)\delta_m^{-n+4}
,&\text{$n$, $m$ odd}.
\end{cases}
\end{equation}
\end{proposition}

\medskip

\begin{remark}
The bilinear form $\gamma(f,g)=\cins fdg$ on the function algebra
$\A=\A_{(e_1,e_2)}$  appearing in the definition of the
cocycle for the current algebra is a cocycle 
for  $\A$ (considered as abelian Lie algebra) defining 
a central extension of it, i.e. an infinite dimensional Heisenberg algebra.
In \cite{Schlloc} it was shown that this cocycle 
is up to multiplication with a constant the
unique  local and multiplicative cocycle for $\A$.
Note that because $\A$ is abelian, there are no non-trivial coboundaries.
A cocycle is called {\it local} if there are constants $M$ and $N$
such that
\begin{equation}\label{E:loc}
\gamma(\A_n,\A_m)\ne 0,\qquad \text{ implies }\quad
M\le n+m\le L.
\end{equation}
Being a local cocycle says, that the central extension  
defined via the cocycle 
is almost-graded.
A cocycle is called {\it multiplicative} if
\begin{equation}\label{E:mult}
\gamma(fg,h)+
\gamma(gh,f)+
\gamma(hf,g)=0,\qquad f,g,h\in\A.
\end{equation}
In the quoted article the uniqueness was shown by applying a recursive
procedure using the locality and the multiplicativity of the
cocycle and the almost-graded structure of the algebra $\A$.
From this it is clear that the total cocycle $\cins fdg$ 
is already fixed by one non-vanishing cocycle value
(let's say $\gamma(A_2,A_{-2})$) and the structure constants 
\refE{structwo} of the algebra.
Instead going through the Laurent series calculations above, one could 
equally well go through the combinatorics of the recursion
in \cite{Schlloc} using the structure of $\A$.
We will supply an alternative  proof in the next appendix. From this proof 
the similarity of the structure equations
\refE{structwo} with the formulas for the central extensions
\refE{famcent1} will become clear.
\end{remark} 

\section{A different proof of \refT{famcent}}\label{S:proofd}
\newcommand{\ga}{\gamma}

First note that as antisymmetric and 
bilinear form, $\ga(f,g)=\cins fdg$  defines 
a two-cocycle for the abelian Lie algebra $\A$.
Furthermore, note that there does not exist any non-trivial coboundary.

Recall that we are working over a compact Riemann surface, resp.
a smooth projective curve over $\C$. By the residue theorem
the total residue of every differential 
has to vanish. Hence, 
a differential having  a pole with a (point-)residue  
has to have at least
another pole.
This implies that  
with respect to the grading introduced, the cocycle is
local (see \refE{loc}). This  can be seen by estimating the pole orders
of the differential $A_ndA_m$ at the points $0$ and $1/2$.

As the (point-)residue of an exact  differential of 
a meromorphic function is zero we obtain
\begin{equation}
0=\cins d(fgh)=\cins fgd(h)+
\cins ghd(f)+
\cins hfd(g).
\end{equation}
This condition is exactly the multiplicativity of the 
cocycle $\ga$.

Recall  the structure of the algebra:
\begin{equation}\label{E:appstruc}
A_n\cdot A_m=
\begin{cases}
A_{n+m},&\text{$n$ or $m$ even,}
\\
A_{n+m}+aA_{n+m-2}+bA_{n+m-4},&\text{$n$ and $m$ odd,}
\end{cases}
\end{equation}
with
$a=3e_1$ and $b=(e_1-e_2)(2e_1+e_2)$.

For the set of  values $\ga(A_n,A_m)$ we introduce its
level $l=n+m$. In the following we will use induction
on the level $l$.

\medskip

{\bf 1.} If $l<0$ then 
$\ga(A_n,A_m)=\ga(A_n,A_{-n+l})=0$.
This follows immediately from the fact that $\ga$ is a local cocycle
(see \cite{Schlloc}, where the inverted grading was used).
Alternatively, it can be easily checked by showing that there
will be no pole at $z=0$ as long as $l<0$.

\medskip

{\bf 2.}
We have to calculate one cocycle value for normalization.
\begin{multline}
\ga(A_2,A_{-2})=- 
\ga(A_{-2},A_{2})=\res_{z=0}
\left((\wp(z)-e_1)^{-1}\frac d{dz}(\wp(z)-e_1)\right)
\\=
\res_{z=0}
\left((\wp(z)-e_1)^{-1}\wp'\right)=-2.\qquad
\end{multline}

\medskip

{\bf 3.}
\begin{lemma} 
\label{L:zero}
$\ga(A_n,A_0)=0$ for all $n\in\Z$.
\end{lemma}
\begin{proof}
\begin{equation}
\ga(A_n,A_0)=\ga(A_n,A_0\cdot A_0)=
-\ga(A_n\cdot A_0,A_0)-\ga(A_0\cdot A_0,A_n)
=-\ga(A_n,A_0)+\ga(A_n,A_0)=0.
\end{equation}
\end{proof}
Above we used the multiplicativity. 
Again with  the multiplicativity 
we obtain
\begin{lemma}\label{L:three}
$\ga(f\cdot f,f)=0$.
\end{lemma}

\medskip

{\bf 4.} Let $l=0$, and $n$ of arbitrary parity.
We calculate
\begin{multline}
\ga(A_{n+1},A_{-(n+1)})=
\ga(A_{n}A_1,A_{-(n+1)})=
-\ga(A_{1}A_{-(n+1)},A_{n})
-\ga(A_{-(n+1)}A_{n},A_{1})
\\
=-\ga(A_{-n},A_{n})
-\ga(A_{-1},A_{1})
=\ga(A_{n},A_{-n})
+\ga(A_{1},A_{-1}).
\end{multline}
Here certain remarks are necessary, as
this chain of equalities is not that innocent as it looks.
 For the first equality note that
if $n$ is even we can indeed replace $A_{n+1}$ by $A_{n}A_1$,
using 
\refE{appstruc}.
If $n$ is odd then  $A_{n+1}=A_{n}A_1-aA_{n-1}-bA_{n-3}$.
But we have $\ga(A_{n-1},A_{-(n+1)})=\ga(A_{n-3},A_{-(n+1)})=0$,
as their levels  are $-2$,  resp. $-4$, and for those levels
the cocycle is vanishing.
Hence, indeed the first inequality is true for every parity.
Using the multiplicativity,  we replace the one cocycle by 
the negative of its two partner cocycles.
Now we use the same kind of arguments as above to replace the 
products for them.

Altogether we obtain a simple recursion relation which has as unique
solution
\begin{equation}\label{lzero}
\ga(A_n,A_{-n})=n\cdot\ga(A_1,A_{-1})=(-n).
\end{equation}
The last equality follows from the cocycle value we calculated in
Step 2.

\medskip

{\bf 5.}
We consider $(n,-n+l)$ with $l\ne 0$, such that either both entries
are even (this says $l$ is even), 
or $n$ is odd and the other entry is even (this says $l$ is odd).
We claim that all cocycle values for such pairs vanish.
For $l$ even, we will show the claim directly. For $l$ odd we
will use induction.
First note that for $l$ negative the values will be zero by locality.
Hence the start of the induction is in the odd case trivially 
true. 
We calculate 
\begin{multline}
\label{E:rec}
\ga(A_{n},A_{-n+l})=
\ga(A_{n-2}A_2,A_{-n+l})=
-\ga(A_{2}A_{-n+l},A_{n-2})
-\ga(A_{-n+l}A_{n-2},A_{2})
\\
=
\ga(A_{n-2},A_{-(n-2)+l})+
\ga(A_{2},A_{l-2}).
\end{multline}
\begin{lemma}\label{L:l2}
For $l\ne 0$ we have 
$ \ga(A_{2},A_{l-2})=0$.
\end{lemma}
\begin{proof}
By the locality $\ga(A_2,A_{l-2})=0$ if $l\le 0$.
First we have to consider $l=1,2,3,4$ individually.
\newline
For $l=1$ we get 
\begin{multline}
\ga(A_2,A_{-1})=
\ga(A_1A_1,A_{-1})=
-\ga(A_{1}A_{-1},A_1)
-\ga(A_{-1}A_{1},A_1)=
\\
-2\ga(A_{-1}A_{1},A_1)=
-2(\ga(A_{0},A_1)
+a\ga(A_{-2},A_1)
+b\ga(A_{-4},A_1))=0.
\end{multline}
Here we used 
that the cocycle vanishes if the level is negative,
the multiplicativity and 
\refL{zero}.
\newline
For $l=2$ we have $\ga(A_2,A_{0})=0$ by \refL{zero}.
\newline
For $l=3$  we calculate
\begin{equation}
\ga(A_2,A_{1})=
\ga(A_1A_1-aA_0-bA_{-2},A_{1})=
\ga(A_1A_1,A_1)=0.
\end{equation}
(We used also \refL{three}.)
\newline
For $l=4$ we have $\ga(A_2,A_{2})=0$ by the antisymmetry.

Now let $l\ge 5$. We  showed above that 
\begin{equation}
\ga(A_{r},A_{-r+l})-
\ga(A_{r-2},A_{-(r-2)+l})+
\ga(A_{2},A_{l-2})=0.
\end{equation}
If $l$ is even we let $r=4,6,\ldots, l-2$ and sum this 
$(l-4)/2$ equations up.
As sum we obtain
\begin{equation}
\ga(A_{l-2},A_{2})- \ga(A_{2},A_{l-2}) -
\frac {l-4}{2} \ga(A_{2},A_{l-2})=0.
\end{equation}
Hence,
\begin{equation} 
\frac {l-4}{2} \ga(A_{2},A_{l-2})=0,
\end{equation}
and as $l>4$, $\ga(A_{2},A_{l-2})=0$.
\newline
For $l$ odd we sum up the equations for $r=3,5,\ldots,l-2$.
These are $(l-3)/2$ equations. The sum calculates to
\begin{equation}\label{E:zw}
\ga(A_1,A_{l-1})+\frac {l-1}{2}\ga(A_2,A_{l-2})=0.
\end{equation}
We calculate
\begin{equation}
\ga(A_2,A_{l-2})=\ga(A_1A_1,A_{l-2})-
a\ga(A_0,A_{l-2})-b\ga(A_{-2},A_{l-2})=
\ga(A_1A_1,A_{l-2}),
\end{equation} as the other 
summands  vanish by
induction.
Hence,
\begin{equation}
\ga(A_2,A_{l-2})=
-\ga(A_1A_{l-2},A_{1})-\ga(A_{l-2}A_1,A_{1})
=2\ga(A_{1},A_{l-1}+aA_{l-3}+bA_{l-5}).
\end{equation}
By induction what remains is
\begin{equation}
\ga(A_2,A_{l-2})=2\ga(A_1,A_{l-1}).
\end{equation}
Hence, from \refE{zw} it follows that also in
this case $\ga(A_2,A_{l-2})=0$.
\end{proof}
This proofs shows also
\begin{lemma}\label{L:l1}
For $l$ odd we have $\ga(A_1,A_{l-1})=0$.
\end{lemma}

If $n$ is even from
\refE{rec} it follows that 
$\ga(A_{n},A_{-n+l})=\frac n2\ga(A_2,A_{l-2})$.
Hence by \refL{l2}, all cocycle values will vanish.

If $n$ is odd, $l$ is odd and we 
get by \refL{l1} and \refL{l2} 
that the starting value and the increment for the recursion are zero.
Hence, also here all cocycle values will vanish.

\medskip

{\bf 6.}
It remains to consider the case
that 
both entries in $(n,-n+l)$ are odd, and $l> 0$.
In any case, $l$ will be an even number.
\begin{multline}
\ga(A_{n},A_{-n+l})=
\ga(A_{n-1}A_{1},A_{-n+l})=
-\ga(A_{1}A_{-n+l},A_{n-1})-
 \ga(A_{-n+l}A_{n-1},A_{1})
\\=
-\ga(A_{-n+l+1}+aA_{-n+l-1}+bA_{-n+l-3},A_{n-1})
-\ga(A_{l-1},A_1).
\end{multline}
In the first expression there are only (even,even) combinations, which
have only a non-vanishing value if they are $(k,-k)$. 
Altogether we obtain 
\begin{equation}
\ga(A_{n},A_{-n+l})=
-a(n-1)\delta_l^2-
-b(n-1)\delta_l^4
+\ga(A_{1},A_{l-1}).
\end{equation}
For $l=2$ we obtain $\ga(A_1,A_{l-1})=\ga(A_1,A_1)=0$.
Hence, $\ga(A_{n},A_{-n+2})=a\cdot(-n+1)$.
\newline
For $l=4$ we calculate 
\begin{equation}
\ga(A_1,A_3)=
\ga(A_1,A_2A_1)=-
\ga(A_1,A_1A_2)-\ga(A_2,A_1A_1).
\end{equation}
This implies
\begin{equation}
\ga(A_1,A_3)=\frac 12\ga(A_1A_1,A_2)=
\frac 12\ga(A_2+aA_0+A_{-2},A_2)=
\frac 12\ga(A_{-2},A_2)=b.
\end{equation}
Hence,  $\ga(A_{n},A_{-n+4})=b\cdot(-n+2)$.
\newline
For $l>4$ we have 
$\ga(A_{n},A_{-n+l})=\ga(A_{1},A_{l-1})$.
This says that for level $l$ all values are equal to the
same constant $\ga(A_{1},A_{l-1})$. 
This is only possible if this constant is zero, as
\begin{equation}
\ga(A_{1},A_{l-1})=- 
\ga(A_{l-1},A_{1})=-
\ga(A_{l-1},A_{-(l-1)+k})=
-\ga(A_{1},A_{l-1}).
\end{equation}
Hence, $\ga(A_{n},A_{-n+l})=0$ for $l>4$.
This closes the calculation of the cocycle values.



\end{document}